


\magnification 1200
\input amstex
\documentstyle{amsppt}
\NoBlackBoxes
\NoRunningHeads
\nologo

\hsize = 6.7 truein
\vsize = 9.3 truein

\redefine\t{^{\sssize T}}

\define\vre{\varepsilon}
\define\hd{Hausdorff dimension}
\define\hs{homogeneous space}
\define\df{\overset\text{def}\to=}
\define\un#1#2{\underset\text{#1}\to#2}
\define\br{\Bbb R}
\define\bn{\Bbb N}
\define\bz{\Bbb Z}
\define\bzp{\Bbb Z_{\sssize +}}
\define\brp{\Bbb R_{\sssize +}}

\define\di{Diophantine}
\define\da{Diophantine approximation}

\define\ve{\bold e}
\define\vu{\bold u}

\define\vx{\bold x}
\define\vy{\bold y}

\define\vv{\bold v}

\define\vs{\bold s}
\define\vw{\bold w}
\define\vp{\bold p}
\define\vq{\bold q}

\define\vf{\bold f}
\define\vt{\bold t}
\define\vg{\bold g}

\redefine\sp{Sprind\v zuk}

\redefine\lfloor{|\langle}
\redefine\rfloor{\rangle|}

\define\nz{\smallsetminus \{0\}}

\define\cag{$(C,\alpha)$-good}

\define\ca{\Cal A(d,k,C,\alpha,\rho)}
\define\sfb{(S,\varphi,B)}
\define\ph{\varphi}
\define\p{\Phi}

\define\ssm{\smallsetminus}

\def\dist{\operatorname{dist}}

\topmatter
\title Khintchine-type theorems 
on manifolds: \\ the convergence case for standard \\ and
multiplicative versions. \endtitle

\author { V.$\,$Bernik, D.$\,$Kleinbock and G.$\,$A.$\,$Margulis\\ {\ }} \endauthor




  

 \thanks The work of the second named author was supported in part by NSF
Grant DMS-9704489, and that of the third named
author by NSF Grant DMS-9800607. \endthanks

    \address{ V.$\,$I.~Bernik, Institute of Mathematics, Belarus
Academy of Sciences, Minsk, 220072, Belarus}
  \endaddress

\email bernik\@im.bas-net.by \endemail

    \address{ D.~Kleinbock,  Department of
Mathematics, Brandeis University, Waltham, MA 02454-9110}
  \endaddress

\email kleinboc\@brandeis.edu \endemail


    \address{ G.$\,$A.~Margulis, Department of Mathematics, Yale University, 
   New Haven, CT 06520}\endaddress

\email margulis\@math.yale.edu
  \endemail

\endtopmatter
\document

\heading{1. Introduction}
\endheading 


The goal of this paper is to prove the convergence part of 
the Khintchine-Groshev Theorem,
as well as its multiplicative version,  for nondegenerate smooth
submanifolds in $\br^n$. The proof combines methods from metric number
theory with a new approach involving the geometry  of lattices in
Euclidean spaces. 
 
\subhead{Notation}\endsubhead
The main objects of this paper are $n$-tuples $\vy = (y_1,\dots,y_n)$
of real numbers viewed as {\sl linear forms\/}, i.e.~as row vectors. In what
follows, $\vy$ will always mean a row vector, and we will be
interested in values of a linear form given by $\vy$ at integer points
$\vq = (q_1,\dots,q_n)\t$, the latter being a column vector. Thus
$\vy\vq$ will stand for $y_1q_1 + \dots + y_nq_n$. Hopefully it will
cause no confusion.

We will study differentiable maps $\vf = (f_1,\dots,f_n)$
from open subsets $U$ of $\br^d$ to $\br^n$; again, $\vf$ will be
interpreted as a row vector, so that $\vf(x)\vq$ stands for $q_1f_1(x)
+ \dots + q_nf_n(x)$. In contrast, the elements of the ``parameter
set'' $U$ will be denoted by $x = (x_1,\dots,x_d)$ without
boldfacing, since the linear structure of the parameter space is not
significant.

For $\vf$ as above we will
denote by $\partial_i \vf:U\mapsto \br^n$, $i = 1,\dots,d$, its partial
derivative (also a row vector) with respect to $x_i$. If $F$ is a
scalar function on $U$, we will denote by $\nabla F$ the column 
vector consisting of  partial derivatives of $F$. With some abuse of
notation, the same way we will treat vector functions $\vf$: namely, $\nabla
\vf$ will stand for 
the matrix function $U\mapsto 
M_{d\times n}(\br)$ with rows given by partial derivatives $\partial_i \vf$. 
We will also need higher 
order differentiation:  for a {\it multiindex\/}
$\beta = (i_1,\dots,i_d)$, $i_j\in \bz_{\sssize +}$, we let $|\beta| =
\sum_{j = 1}^d i_j$ and $\partial_\beta =
\bold \partial_1^{i_1}\circ\dots\circ \partial_d^{i_d}$. 


Unless otherwise indicated, the norm 
$\|\vx\|$ of a vector $\vx\in\br^k$ (either row or column vector) 
will stand for $\|\vx\| = \max_{1\le i \le k}|x_i|$. In some cases
however we will work with the Euclidean norm $\|\vx\| = \|\vx\|_e =
\sqrt{\sum_{i 
= 1}^k x_i^2}$, 
keeping the same notation. This distinction 
will be clearly emphasized to avoid confusion. 
We will denote by $\br^k_1$ the set of unit vectors in $\br^k$ (with
respect to the Euclidean norm).

We will use the notation $\lfloor x
\rfloor$ for 
the distance between $x\in\br$ and the closest integer, $\lfloor x \rfloor
\df \min_{k\in\bz} |x - k|$. (It is quite customary to use $\|x\|$
instead, but we are not going to do this in order to save the latter
notation for norms in vector spaces.)  If $B\subset
\br^k$, 
we let $|B|$ stand for the Lebesgue
measure of $B$.

\subhead{Basics on \da}\endsubhead In what follows, we let 
$\Psi$ be a positive function   defined on $\bz^n\nz$, and consider the set
$$
\Cal W(\Psi) \df \big\{\vy\in\br^n\bigm| \lfloor \vy\vq\rfloor \le \Psi(\vq)\text{ for infinitely many }\vq\big\}\,.
$$
Clearly the faster $\Psi$ decays at infinity, the smaller is the set $\Cal W(\Psi)$. In particular, the Borel-Cantelli Lemma gives a sufficient condition for this set to have measure zero: $|
\Cal W(\Psi)| = 0$ if 
$$
\sum_{\vq\in \bz^n\nz} {\Psi(\vq)} < \infty\,.\tag 1.1
$$
It is customary to refer to the above statement, as well as to its various analogues, as to the convergence case of a Khintchine-type theorem, since 
the fact that (under some regularity restrictions on $\Psi$) the condition (1.1) is also necessary was first proved by A.~Khintchine for $n = 1$ and later generalized by A.~Groshev  and W.~Schmidt. See \S 8.5 for more details. 

The standard class of  examples 
is given by functions which depend only on the norm of $\vq$. 
If $\psi$ is a positive function defined on positive integers, let us say, following \cite{KM2}, that $\vy$    is {\sl
$\psi$-approximable\/},  to be abbreviated as {\sl $\psi$-A\/}, if
it belongs to $\Cal W(\Psi)$ where  \footnote{We are grateful to M.$\,$M.$\,$Dodson for
 a permission to modify his terminology used in
\cite{Do} and \cite{BD}. In our opinion, the parametrization (1.2s) instead of the traditional $\psi(\|\vq\|)$ makes the structure  
more transparent and less dimension-dependent; see
\cite{KM1, KM2} for justification.}
$$
\Psi(\vq) =\psi(\|\vq\|^{n})\,.\tag 1.2s
$$ 
If 	$\psi$ is non-increasing (which will be our standing assumption), (1.1) is satisfied 
 if and only if  
$$
\sum_{k = 1}^{\infty}{\psi(k)} < \infty\,.\tag 1.1s
$$
 An example: 
almost all   $\vy\in\br^n$ are not VWA ({\sl very well
approximable\/}, see \cite{S2, Chapter IV, \S 5});
the latter is defined to be $\psi_\vre$-approximable for some
positive $\vre$, with 
$\psi_\vre(k)\df k^{-(1+\vre) }$.

Another important special case is given by 
$$
\Psi(\vq) =\psi\big(\Pi_{\sssize +}(\vq)\big)\,.\tag 1.2m
$$ 
where
$
\Pi_{\sssize +}(\vq)$ is defined as  $\prod_{i = 1}^n \max(|q_i|, 1)$,
 i.e.~the 
absolute value of the product of all the nonzero coordinates of
$\vq$. We will say  that $\vy\in\br^n$ is {\sl $\psi$-MA\/} ({\sl $\psi$-multiplicatively approximable})   if
it belongs to $\Cal W(\Psi)$ with $\Psi$ as in (1.2m). In this case, again assuming the monotonicity of $\psi$,  (1.1) is satisfied 
if and only if 
$$\sum_{k = 1}^{\infty}(\log k)^{n-1}{\psi(k)} < \infty\,.\tag 1.1m
$$ Also, 
since $\Pi_{\sssize +}(\vq)$ is not greater than $\|\vq\|^n$,
 any
$\psi$-approximable $\vy$ is automatically 
$\psi$-MA.  For example, one can define {\sl very well
multiplicatively approximable\/} ({\sl VWMA})  points to be
$\psi_\vre$-multiplicatively approximable for some 
positive $\vre$, with $\psi_\vre$ as above; it follows that almost all
$\vy\in\br^n$ are not VWMA.

\subhead{\da\ on manifolds}\endsubhead Much more  intricate questions arise if one restricts $\vy$ to lie  on
a submanifold $M $ of $\br^n$. In 1932 K.~Mahler \cite{M}
conjectured that almost all points on the curve  
 $$
\{(x,x^2,\dots,x^n)\mid x\in \br\}\tag 1.3
$$ 
are not VWA. V.~Sprind\v zuk's  proof of this conjecture (see
\cite{Sp1, Sp2}) has eventually led to the development of a new branch
of metric number theory, usually referred to as ``\da\ with dependent
quantities'' or ``\da\ on manifolds''. In particular, Sprind\v zuk's
result was improved by A.~Baker \cite{B1} in 1966: he showed that if
$\psi$ is a 
positive  non-increasing
function such that   
$$
\sum_{k = 1}^\infty \frac{\psi(k)^{1/n}}{k^{1 - 1/n}} < \infty\,,\tag 1.4
$$
then  almost all points on the
curve  (1.3) 
are not $\psi$-approximable.  Baker conjectured that (1.4) could be replaced by
the optimal condition (1.1s); this conjecture  was proved later by
V.~Bernik \cite{Bern}. 
As for the multiplicative approximation, it was conjectured by
A.~Baker in his book \cite{B2} that   almost all points on the
curve  (1.3) are not VWMA; the validity of this conjecture for   $n
\le 4$ was verified in 1997 by V.~Bernik and V.~Borbat  \cite{BB}. 

Since the mid-sixties, a lot of efforts have been directed to obtaining
similar results for larger classes of smooth submanifolds of
$\br^n$. A new method, based on combinatorics of the space of
lattices, was developed  in 1998 in the paper
\cite{KM1} by  Kleinbock and Margulis. Let us employ the following
definition: if $U$ is an open subset of
 $\br^d$ and $l\le m \in \bn$, say that an $n$-tuple  $\vf =
(f_1,\dots,f_n)$ of  $C^m$ 
functions $U\mapsto \br$ is {\sl $l$-nondegenerate at $x\in U$\/} if 
the space $\br^n$ is spanned by partial derivatives of $\vf$ at
$x$ of order up to $l$.  We will say that $\vf$ is 
 {\sl nondegenerate\/} at $x$ if it is $l$-nondegenerate for some
$l$. 
If $M\subset \br^n$ is a $d$-dimensional $C^m$ submanifold, we will say
that $M$ is {\sl nondegenerate at $\vy\in M$} if any (equivalently,
some)  diffeomorphism
$\vf$ between an open subset $U$  of $\br^d$ and a neighborhood of
$\vy$ in $M$ is   nondegenerate at $\vf^{-1}(\vy)$. We will say that
$\vf:U\to \br^n$ (resp.~$M\subset \br^n$) is {\sl nondegenerate\/} 
if it is nondegenerate at almost every point of $U$  (resp.~$M$, in
the sense of the 
natural measure class on $M$). 

\proclaim{Theorem A {\rm \cite{KM1}}} Let $M$ be a nondegenerate
$C^m$ submanifold of $\br^n$. Then almost all points of $M$  
	are not VWMA (hence not VWA as well).
\endproclaim

In particular, the aforementioned multiplicative conjecture of Baker
follows from this theorem.  
Note also that if the functions $f_1,\dots,f_n$ are analytic and $U$
is connected, 
the nondegeneracy of $\vf$ is equivalent to  the linear independence of 
$1,f_1,\dots,f_n$ over $\br$; in this setting the above statement was
conjectured by \sp\ \cite{Sp4, Conjectures H$_1$, H$_2$}.  


\subhead{Main results and the structure of the paper}\endsubhead The
primary  goal of the 
present paper is to obtain a Khintchine-type 
generalization of Theorem A. 
More precisely,
we prove the following

\proclaim{Theorem 1.1} Let  $U\subset \br^d$ be an open set and
let $\vf:U\to \br^n$  be a  nondegenerate $n$-tuple of $C^m$ functions
on $U$.  
Also let $\Psi: \bz^n\nz \mapsto(0,\infty)$ be a function satisfying {\rm (1.1)} and such that for $i = 1,\dots,n$ one has  
$$
\Psi(q_1,\dots,q_i,\dots,q_n) \ge \Psi(q_1,\dots,q'_i,\dots,q_n) \quad\text{whenever} \quad |q_i| \le |q'_i|\text{ and }q_iq'_i > 0\tag 1.5
$$
(i.e., $\Psi$ is non-increasing with respect to the absolute value of any
coordinate in any orthant of $\br^n$). Then 
$
|\{x\in U\mid \vf(x)\in\Cal W(\Psi)\}| = 0$. 
\endproclaim

In particular, if $\Psi$ is of the form (1.2s) or (1.2m) for a non-increasing function $\psi:\bn\mapsto \br_+$, condition (1.5) is clearly satisfied. Thus 
one has

\proclaim{Corollary 1.2} Let  
$\vf:U\to \br^n$  be as in Theorem 1.1, 
and let $\psi:\bn\mapsto(0,\infty)$ be a
non-increasing function.  Then:

\roster
\item"(S)" assuming {\rm (1.1s)},   
  for        
almost all $x\in U$ the points $\vf(x)$ 
are not $\psi$-A;

\item"(M)" assuming {\rm (1.1m)},   
  for
almost all $x\in U$ the points $\vf(x)$ 
are not $\psi$-MA.

\endroster
\endproclaim

It is worth mentioning that the statement (S)  was recently proved in a paper
\cite{Be5} of V.~Beresnevich using a refinement of \sp's method of ``essential and inessential domains''. Earlier several  special cases were  treated in \cite{DRV1, BDD, Be2}. A preliminary version \cite{BKM} of the present paper, where the two statements of Corollary 1.2  were proved for the case $d = 1$, appeared in 1999 as a
preprint of the University of Bielefeld. 

%

Our proof  of Theorem 1.1 
is based on carefully
measuring sets of 
solutions of certain systems of \di\ inequalities. Specifically, we  fix
a ball $B\subset \br^d$  and look at the set of all $x\in B$ for which
there exists an integer vector $\vq$ in a certain range such that the
value of the function $F(x) =
\vf(x)\vq$ is close to an integer. Our estimates will require
considering two special cases: when the norm of the
gradient $\nabla F(x) = \nabla\vf(x)\vq$
is big, or
respectively, not very big. We will show  in
\S 8.1 that, by means of  
straightforward measure 
computations, 
Theorem 1.1 reduces to 
the following two theorems:

\proclaim{Theorem 1.3} Let $B \subset\br^d$ be a ball of radius $r$,
let  $\tilde B$ stand for the ball with the same center as
$B$ and  of radius $2r$, and let 
functions 
$\vf= (f_1,\dots,f_n)\in C^2(\tilde B)$ be given. Fix $\delta > 0$ and define
$$
L = 
\max_{|\beta| = 2,\,x\in \tilde B}\|\partial_\beta
\vf(x)\|\,.\tag 1.6a 
$$
 Then for every $\vq\in\bz^n$
such that 
$$
\|\vq\| \ge \frac1{4nLr^2}\,,\tag 1.6b
$$
 the set of solutions $x\in B$ of the
inequalities  
$$
\lfloor \vf(x)\vq \rfloor  < \delta \tag 1.6c
$$
and 
$$
 \|\nabla\vf(x)\vq\| \ge \sqrt{ndL\|\vq\|}
\tag 1.6d
$$
has measure at most $C_d\delta|B|$, where $C_d$ is a constant dependent
only on $d$. 
  \endproclaim

\proclaim{Theorem 1.4} Let $U\subset\br^d$ be an open set, $x_0\in U$, and
let $\vf = (f_1,\dots,f_n)$
be an   $n$-tuple of smooth functions on $U$ which is $l$-nondegenerate at
$x_0$. Then there exists a neighborhood $V\subset U$ of  $x_0$ with the following property: for
any ball $B\subset V$ there exist
$E > 0$ such that for any choice of  
$$
0 < \delta \le 1, \quad {T}_1, \dots, {T}_{n}\ge 1\quad \text{and} \quad  K > 0 \text{ with }\frac{\delta K {T}_1\cdot\dots\cdot
{T}_{n}}{\max_i{T}_i} \le 1\,,\tag 1.7a 
$$
 the set 
$$
\Big\{x\in B\bigm| \exists\,\vq\in\bz^n\nz\text{ such that }\cases 
\lfloor \vf(x)\vq \rfloor  < \delta \\
 \|\nabla\vf(x)\vq\|  < K\\
 |q_i| < {T}_i,\quad i = 1,\dots,n
 \endcases\ \Big\}\tag 1.7b
$$
has measure at most $E \vre^{\frac1{d({2l-1})}}|B|$, where one
defines
$$
\vre \df \max\left(\delta, \left(\frac{\delta K {T}_1\cdot\dots\cdot
{T}_{n}}{\max_i{T}_i}\right)^{\frac1{n+1}}\right)\,.\tag 1.7c
$$
  \endproclaim

Theorem 1.3, roughly speaking, says that  a function
with  big gradient and not very big second-order partial  derivatives cannot
have values very close to integers on a set of big measure. It  is
proved in \S 2 using an argument  which is 
apparently originally due to Bernik, and is (in the case $d = 1$)
implicitly contained in 
one of the steps of the paper \cite{BDD}. 

As for Theorem 1.4, it is done by a modification of a method  from
\cite{KM1} involving  the geometry of lattices in Euclidean spaces.
The connection with lattices is discussed in \S 5, where  Theorem 1.4
is translated into the language of lattices (see Theorem 5.1). To prove the
latter we
 rely on the notion of \cag\ functions introduced in
\cite{KM1}. This concept is reviewed in \S 3, where, as well as in \S
4, we prove that 
certain functions arising in the proof of 
Theorem 5.1 are \cag\ for suitable $C,\alpha$.  We prove Theorem 5.1
in \S 7, after doing some preparatory work in the preceding section. The
last section of the paper is devoted to the reduction of Theorem 1.1 to Theorems 1.3 and 1.4, as well as to 
several concluding remarks, 
applications  and
some open questions. In particular, there we discuss the complementary 
divergence case of Theorem 1.1, and also present applications involving 
approximation of zero by values of functions and their derivatives.

%

\heading{2. Proof of Theorem 1.3}
\endheading 

We first state a covering result which is well known (and is one of
the ingredients of the main estimate from \cite{KM1}).

\proclaim{Theorem 2.1 (Besicovitch's Covering Theorem{\rm,  see
\cite{Mat, Theorem 2.7}})}   There is an integer $N_d$ depending only
on $d$ with the following property: let  $S$ be a bounded subset of
$\br^d$ and let $\Cal B$ be a family  of nonempty open balls in
$\br^d$ such that each $x\in S$ is the center of some ball of $\Cal
B$;  
then there exists a finite or countable subfamily $\{U_i\}$ of $\Cal B$ with $
1_S \le \sum_i 1_{U_i} \le N_{d}
$ (i.e., $S\subset \bigcup_i U_i$ and the multiplicity of that
subcovering is at most $N_{d}$). 
\endproclaim

This theorem, and the constant $N_d$, will be used repeatedly in the
paper. 

\proclaim{Lemma 2.2} Let   $B \subset\br^d$ be a ball of radius $r$,
and let the 
numbers 
$$
M \ge 1/{4r^2}\tag 2.1a
$$
and $\delta > 0$ be given. 
Denote by $\tilde B$ the ball with the same center as
$B$ and  of radius $2r$. 
Take a function 
$F\in C^2(\tilde B)$ such that 
$$
\sup_{|\beta| = 2,\,x\in \tilde B}|\partial_\beta F(x)| \le
{M}\,, \tag 2.1b
$$ 
and denote by $S$ the set of all $x\in B$ for which the inequalities
$$
\lfloor F(x) \rfloor  < \delta \tag 2.1c 
$$
and
$$
\|\nabla F(x)\| \ge \sqrt {dM}\tag 2.1d 
$$
hold. Then $|S| \le C_d\delta |B|$, where $C_d$ is a constant dependent
only on $d$. \endproclaim

\demo{Proof} Clearly $|S| \le 16\delta |B|$ when $\delta \ge
1/16$, so without loss of generality we can assume that 
$\delta$ is less than $1/16$. 
Also, given $x\in S$, without loss of generality we can assume that
the maximal value of $|\partial_j F(x)|$, $j =
1,\dots,d$, occurs when $j = 1$. Denote
$\frac1{2|\partial_1 F(x)|}$  by $\rho$; note that
$\rho\sqrt d \le \frac{1}{ 2\sqrt{M}} \le r$ due to (2.1a) and
(2.1d), therefore 
the ball $B(x, \rho\sqrt d)$  is contained in $\tilde
B$. Also let us denote by 
$U(x)$ the maximal 
ball centered in $x$  such that $\lfloor
F(y) \rfloor < 
1/4$ for all $y\in U(x)$. It is clear that  there exists a unique $p \in\bz$
such that 
$
| F(y) + p |  < 1/4
$  for all $y\in U(x)$. 
We claim  that the radius of $U(x)$ is not bigger than
$\rho$. Indeed, 
one has
$$
F(x_1\pm\rho,x_2,\dots,x_d) + p = F(x) + p \pm
{\partial_1 F(x)}\rho + \tfrac{\partial^2_1F(z)}2\rho^2
$$
for some $z$ between $x$ and $(x_1\pm\rho,x_2,\dots,x_d)$.
Thus 
$$
|F(x_1\pm\rho,x_2,\dots,x_d) + p| \un{(2.1bc)}{\ge} -\delta + 1/2 -
M\rho^2/2
\un{(2.1d)}\ge 1/2 - \tfrac{1}{8\sqrt{d}} - \delta \un{($\delta < 1/16$)}{>} 1/4\,,
$$
 and the claim is proved. In particular, $U(x)\subset \tilde B$; moreover,
if one denotes by $\bar U(x)$ the cube circumscribed around $U(x)$ with
sides parallel to the coordinate axes, then $\bar U(x)\subset \tilde
B$ as well. 

On the other hand, the radius of $U(x)$ cannot be too small: if
$y\in B(x, \tfrac\rho{4\sqrt d})$, one has
$$
|F(y) + p| \le | F(x) + p | + \left|\tfrac{\partial F}{\partial
u}(x)\right|\tfrac{\rho}{4\sqrt d} + \tfrac12 \tfrac{\partial^2 F}{\partial
u^2}(z)\tfrac{\rho^2}{16d}\,,
$$
where $u$ is the unit vector parallel to $y - x$, and $z$ is
between  $x$ and $y$. 
Note that it follows from our ordering of coordinates that
$\left|\tfrac{\partial F}{\partial 
u}(x)\right| \le \sqrt d |\partial_1 F(x)| = \tfrac{\sqrt d}{2\rho}$, and from (2.1b)
that $\left|\tfrac{\partial^2 F}{\partial
u^2}(z)\right| \le d M \le \tfrac{d}{4\rho^2}$. Therefore
$
|F(y) + p| \le  \delta + 1/8 +  
1/128 < 1/4\,,
$
which shows that $U(x)\subset B(x, \tfrac\rho{4\sqrt d})$, and, in
particular, 
$$
|U(x)| \ge 
C'_d \rho^d\tag 2.2
$$ 
(the values of constants $C'_d$, and also of $C''_d$, $C'''_d$ to be
introduced later, 
depend
only on $d$). 

Also one can observe that 
$\partial_1 F(y)$ does not oscillate too much when $y\in \bar
U(x)$:  for  
some $z$ between  $x$ and $y$
one gets 
$$
|\partial_1 F(y) - \partial_1 F(x)| \le  |\tfrac{\partial}{\partial
u}\partial_1 F (z)|\rho\sqrt d \un{(2.1b)}{\le} {M}\rho\sqrt d
\un{(2.1d)}{\le} 1/4\rho = 
 \tfrac{|\partial_1 F(x)|}2\,,
$$
(here again $u$ is the unit vector parallel to $y - x$). This
implies that  the absolute value of ${\bold \partial}_1 F(y)$, $y\in \bar U(x)$,
is not less than $\tfrac12|\partial_1 F(x)|$; in
particular, for every $y_2,\dots,y_d$ such that $|y_i - x_i | < \rho$,
$i = 2,\dots,d$, the function $F(\cdot,y_2,\dots,y_d)$  is monotonic on $(x_1 -
\rho,x_1 + \rho)$, and therefore
$$
\big|\{y_1 \in (x_1 -
\rho,x_1 + \rho)\bigm| |F(y_1,\dots,y_d) + p| < \delta\}\big| \le  {2\delta}
\tfrac{2}{|{\partial}_1 F(x)|} = {8\rho\delta}\,.
$$
Now we can estimate $\big|\big\{x\in U(x)\bigm| | F(x) + p |  <
\delta \big\}\big|$ from above by
$$
\big|\big\{x\in \bar U(x)\bigm| | F(x) + p |  < \delta \big\}\big| \le
8\rho\delta 2^{d-1}\rho^{d-1} = C''_d \delta \rho^d
\un{(2.2)}{\le}C'''_d|U(x)|\,. 
$$
The set $S$ is covered by all the balls $U(x)$, $x\in S$, and, using
Theorem 2.1,  one can choose a subcovering $\{U_i\}$ of
multiplicity  at most $N_{d}$. Then one has
$$
|S| \le \sum_i C'''_d \delta |U_i|\le C'''_d N_d \delta|\tilde B|  = C_d \delta|B|
\,,
$$  
which finishes the proof. 
\qed
\enddemo

Now it takes very little to  complete the

\demo{Proof of Theorem 1.3} Given the balls $B\subset \tilde B\subset
\br^d$, an 
$n$-tuple of $C^2$ functions $\vf$ on $\tilde B$,  a
positive $\delta$ and 
 $\vq\in\bz^n$ satisfying
(1.6b) with $L$ as in (1.6a), denote $
F(x)\df \vf(x)\vq$ and $M \df nL\|\vq\|$. Then inequalities
(1.6abcd) can be rewritten as (2.1bacd), and the theorem follows.  
 \qed\enddemo

\heading{3. \cag\ functions}
\endheading 
Let us recall  the  definition introduced in
\cite{KM1}. If  $C$ and $\alpha$ are  positive
numbers and $V$ a subset of $\br^d$, let us say that   a function
$f:V\mapsto \br$ is {\sl \cag\ on\/}   $V$  if for any open ball
$B\subset V$ and any $\vre > 0$
one has 
$$
\big|\{x\in B\bigm| |f(x)| < \vre\cdot{\sup_{x\in
B}|f(x)|}\}\big| \le 
C\vre^\alpha |B|\,. 
\tag 3.1 
$$

Several elementary facts about \cag\ functions are listed below:

 \proclaim{Lemma 3.1} {\rm (a)} $f$ is \cag\ on $V$ $\Rightarrow$  so
is $\lambda f$ $\forall\,\lambda\in \br$; 

{\rm (b)}  $f_i$, $i\in I$, are \cag\ on $V$
$\Rightarrow$ 
so is $\sup_{i\in I}|f_i|$;

{\rm (c)} If $f$ is  \cag\ on $V$ and $c_1\le \tfrac
{|f(x)|}{|g(x)|}\le c_2$
for all $x\in V$, then $g$ is \newline
$\big(C(c_2/c_1)^\alpha,\alpha)$-good on $V$; 

{\rm (d)} $f$ is  \cag\ on $V$ $\Rightarrow$ it is $(C',\alpha')$-good
on $V'$ for every $C'\ge C$, $\alpha' \le \alpha$ and $V'\subset V$. 
 \endproclaim  

Note that it follows from part (b) that the (supremum) norm of a
vector-function $\vf$ is \cag\ whenever every component of $\vf$ is
\cag. Also part (c) shows that one is allowed to replace the norm by
an equivalent one, only affecting 
$C$ but not $\alpha$. 

 \proclaim{Lemma 3.2} Any polynomial $f\in\br[x_1,\dots,x_d]$ of
degree not greater than ${l}$ is  
$(C_{d,{l}}, \frac1{d{l}})$-good on $\br^d$, where $C_{d,{l}} =
\frac{2^{d+1}d{l}({l}+1)^{1/{l}}}{v_d}$ (here and in the next lemma  $v_d$ stands for the volume
of the unit 
ball in $\br^d$). 
 \endproclaim  

\demo{Proof} The case $d = 1$ is proved in \cite{KM1, Proposition
3.2}. By induction on $d$, as in the proof of \cite{KM1, Lemma
3.3}, one can show that  for any $d$-dimensional cube $B$ and for any
$\vre > 0$ one has 
$$
\big|\big\{x\in B\bigm| |f(x)| < \vre\cdot \sup_{x\in
B}|f(x)|\big\}\big| \le 2d{l}({l}+1)^{1/{l}} \vre^{1/d{l}}|B|\,,
$$
and the claim follows by circumscribing a cube around any ball in
$\br^d$. \qed\enddemo

The next lemma is a direct consequence of \cite{KM1,
Lemma 3.3}.

\proclaim{Lemma 3.3}   Let $U$ be an open subset of $\br^d$,  and let $f\in C^k(V)$ be such that for some constants $A_1, A_2 > 0$ one has 
$$
|\partial_\beta f(x)| \le A_1 \quad \forall\,\beta\text{ with }|\beta|\le k\,,\tag 3.2$\le$
$$
and 
$$
|\partial_i^{k}f(x)|\ge A_2\quad \forall\, i = 1,\dots,d\tag 3.2$\ge$
$$
for all $x\in U$. Also  let $V$ be a subset of $U$ such that whenever a ball $B$ lies in $V$, any cube 
circumscribed around $B$ is contained in $U$. Then $f$ is  $(C,\frac1{d{l}})$-good   on $V$, where \newline $C 
= \frac{2^d}{v_d}dk(k+1)\left(\frac{A_1}{A_2}(k+1)(2k^k+1)\right)^{1/k}$. 

\endproclaim  

The following proposition describes the \cag\
property of functions chosen from certain compact families defined by
non-vanishing of partial derivatives. The argument used is similar to
that of Proposition 3.4 from \cite{KM1}.

\proclaim{Proposition  3.4}   
Let $U$ be an open subset of $\br^d$,  and let 
$\Cal F\subset C^{l}(U)$ be a  family of functions 
$f:U\mapsto \br$ such that 
$$
\{\nabla f\mid f\in \Cal F\}\text{ is compact in }C^{{l}-1}(U)\,.\tag 3.3
$$
 Assume also that 
$$
\inf \Sb{f\in \Cal F } \endSb 
\sup_{|\beta|\le {l}}
|\partial_\beta f(x_0)| > 0\tag 3.4
$$
(in other words, the derivatives of $f$ at $x_0$ uniformly 
(in $f\in \Cal F$) generate $\br$).
Then there exists  a
neighborhood $V\subset U$ of $x_0$ and a positive $C = C(\Cal F)$
such that the following holds for all $f\in \Cal F$:

{\rm (a)}  $f$ is  $(C,\frac1{d{l}})$-good   on $V$;

{\rm (b)}  $\|\nabla f\|$ is  $(C,\frac1{d({l}-1)})$-good   on $V$.
\endproclaim

\demo{Proof} Assumption (3.4) says that there exists a constant $C_1 >
0$ such that for any $f\in \Cal F$ one can find a multiindex $\beta$
with $|\beta| = k \le {l}$ with 
$$
\big|\partial_\beta f(x_0)\big| \ge C_1\,.
$$
 By an appropriate rotation  of the coordinate system 
one can guarantee that  
$|\partial_i^k f(x_0)|\ge C_2$ for  all $i = 1,\dots,d$ and some
positive $C_2$ independent of $f$. Then one uses the continuity of the
derivatives of $f$ and compactness of $\Cal F$ to choose a
neighborhood $V'\subset U$ of $x_0$ 
and positive $A_1$, $A_2$  (again independently of $f\in \Cal F$) such that the inequalities (3.2) hold 
for all $x\in V'$ (note that in the above
inequalities  both $k\in\{1,\dots,l\}$ and the coordinate system
depend on $f$). Finally we let $V$ be a smaller neighborhood of 
$x_0$ such that whenever a ball $B$ lies in $V$, any cube $\hat B$
circumscribed around $B$ is contained in $V'$. 

Now part (a) immediately follows from the previous lemma. As for the second part, let $\Cal F'$ be the family of
functions in $\Cal F$ such that  (3.2$\ge$) holds with $k \ge 2$. Then
the family $\{\partial_i f\mid f \in\Cal F',\,i = 1,\dots,d\}$ 
satisfies (3.3)  and (3.4) with $l-1$ in place of $l$; hence, by part
(a), all functions from this family are $(C',\frac1{d({l}-1)})$-good on
some neighborhood of $x_0$ with
some uniform constant $C'$. Thus, by virtue of Lemma 3.1(b,c), the norm
of $\nabla f$ is $(C,\frac1{d({l}-1)})$-good with perhaps a different constant $C$. 

It remains to consider the case when $k$ as in
(3.2$\ge$) is equal to $1$. Then $A_1 \le \|\nabla f(x)\| \le A_2$ 
for $x\in B$, therefore 
for any
positive $\vre$ and any $B\subset V$ one has 
$$
 \big|\{x\in B\bigm| \|\nabla f(x)\| < \vre\cdot\sup_{x\in
B}\|\nabla f(x)\|\}\big| \le 
\left(\frac{A_1}{A_2}\right)^{\frac1{n-1}}\vre^{\frac1{n-1}}  |B|\,. \quad\qed
$$
\enddemo

\proclaim{Corollary 3.5}  Let $U$ be an open
subset of $\br^d$, $x_0\in 
U$, and let $\vf = (f_1,\dots,f_n):U\mapsto \br^n$ 
 be an  $n$-tuple of smooth functions which is 
$l$-nondegenerate at $x_0$. Then there exists  a
neighborhood $V\subset U$ of $x_0$ and a positive $C$ such that

{\rm (a)} any
linear combination of $ 1,f_1,\dots,f_n$ is $(C,\frac1{dl})$-good on $V$;

{\rm (b)} the norm of  any
linear combination of $\nabla f_1,\dots,\nabla f_n$ is
$(C,\frac1{n-1})$-good on $V$. 
\endproclaim

\demo{Proof} Take $f = c_0 + \sum_{i = 1}^nc_if_i$; in view of Lemma
3.1(a), one can without loss of 
generality assume that the norm of $(c_1,\dots,c_n)$ is equal to
$1$. All such functions $f$ belong to a family satisfying
(3.3) and (due to the nondegeneracy of $\vf$ at $x_0$) (3.4), thus
the above proposition applies.
\qed\enddemo

We close the section with two auxiliary lemmas which will be used below
to prove that certain functions are \cag. 

\proclaim{Lemma 3.6}  Let $B$ be a ball in $\br^d$ of radius $r$, and let 
$f \in C^{l}(B)$ and $c > 0$ be  such that for some unit vector $u$
in $\br^d$, some $k\le l$   and all $x\in B$ one has 
$\big|\frac{\partial^k f}{\partial u^k}(x)\big| \ge c$.
Then   
$$
\sup_{x,y\in B}|f(x) - f(y)| \ge
\frac{c}{k^k(k+1)!}(2r)^k\,.\tag 3.5 
$$
\endproclaim

\demo{Proof} If $x_0$ is the center of $B$, consider the function
$g(t) = f(x_0 + tu)$ defined on $I \df [-r,r]$.  
Denote $\sup_{s,t\in I}|f(s) - f(t)|$ by $\sigma$. We claim that
$$
\sigma \ge
\frac{a}{k^k(k+1)!}(2r)^k\,;\tag 3.6 
$$
this clearly implies (3.5). To prove (3.6), take any $s\in I$, divide
$I$ into  $k$ equal segments  and let $p(t)$ be the Lagrange
polynomial of degree $k$ formed by using values of $g(s) - g(t)$ at
the boundary points of these segments. Then there exists $t\in I$ such
that $p^{(k)}(t) = g^{(k)}(t)$, hence, by the assumption,
$|p^{(k)}(t)| \ge c$. On the other hand, after differentiating $p(t)$
$k$ times (see \cite{KM1, (3.3a)}) one gets   
$|p^{(k)}(t)| \le (k + 1) \dfrac{\sigma k!}{(2r/k)^k}$. Combining the
last two inequalities, one easily gets the desired estimate. 
 \enddemo

\proclaim{Lemma 3.7} Let $V\subset \br^d$ be an open ball, and let
$\tilde V$ be the ball with the same center as 
$V$ and twice bigger radius. Let $f$ be a continuous function on
$\tilde V$, and suppose $C,\alpha > 0$ and $0 < \delta < 1$ are
such that  
{\rm (3.1)} holds for any ball $B\subset \tilde V$ and any $\vre \ge
\delta$. Then $f$ is $(C,\alpha')$-good on $V$ whenever $0<\alpha' <
\alpha$ is such that 
$CN_d\delta^{\alpha - \alpha'} \le 1$ (here $N_d$ is the constant
from Theorem 2.1).    
\endproclaim

\demo{Proof} Take  $\vre > 0$ and a ball $B\subset V$,  and denote 
$$
S_{B,\vre}\df\{x\in \tilde V\bigm| |f(x)| < \vre\cdot{\sup_{x\in
B}|f(x)|}\}\,.
$$ The goal is to prove that the measure of $B\cap
S_{B,\vre}$ is not greater than $C\vre^{\alpha'}  |B|$. 

Obviously it suffices to consider $\vre < 1$. Choose $m\in \bzp$ such that
$\delta^{m+1} \le  
\vre < \delta^m$. We will show by induction on $m$ that
$$
|B\cap
S_{B,\vre}| \le
C^{m+1}N_d^m\vre^\alpha  |B|\,. 
\tag 3.7
$$
Indeed, the case $m= 0$ follows from the assumption. Assume that (3.7) holds 
for  
some $m$, and  for every $y\in B\cap
S_{B,\vre}$ let $B(y)$ be the maximal ball
centered in $y$ 
and contained in $S_{B,\vre}$. Observe that, by the continuity of $f$, one has
 $\sup_{x\in B(y)}|f(x)| = \vre\cdot{\sup_{x\in B}|f(x)|}$
for every $y\in B$. Clearly the set $B\cap
S_{B,\vre}$ is covered by all the balls
$B(y)$, and, using 
Theorem 2.1,  one can choose a subcovering $\{B_i\}$ of
multiplicity  at most $N_{d}$. Therefore one has
$$
\split
|B\cap
S_{B,\delta\vre}| 
&\le \sum_i \big|\{x\in B_i\bigm| |f(x)| < \delta\sup_{x\in B_i}|f(x)|\}\big| 
\le \sum_i C\delta^\alpha |B_i| \\
&= 
C\delta^\alpha N_d |B\cap
S_{B,\vre}| 
\le CN_d \cdot 
C^{m+1}N_d^m(\delta\vre)^\alpha |B| 
\,,
\endsplit
$$
which proves (3.7) with $\delta\vre$ in place of $\vre$.

It remains to  write
$$
C^{m+1}N_d^m\vre^\alpha 
= C\cdot (CN_d)^m (\vre)^{\alpha - \alpha'}   \vre^{\alpha'}
< C\cdot (CN_d)^m (\delta^m)^{\alpha - \alpha'} \vre^{\alpha'}
= C\cdot (C N_d\delta^{\alpha - \alpha'})^m
\vre^{\alpha'}
\,,
$$
which implies that $f$ is $(C,\alpha')$-good on $B$ provided
$CN_d\delta^{\alpha - \alpha'} \le 1$. 
\qed\enddemo


\heading{4. Skew-gradients}
\endheading 
In this section we will define and study the following
construction. The main object will be a pair of real-valued differentiable functions
$g_1,g_2$ defined on an open subset $V$ of 
$\br^d$, that  is, a map $\vg:V\mapsto \br^2$. For such a pair,
let us define its {\sl skew-gradient\/} $\tilde\nabla \vg: V\mapsto
\br^d$ by
$$
\tilde\nabla \vg(x)\df g_1(x) \nabla g_2(x) - g_2(x)
\nabla g_1(x)\,.
$$
Equivalently, the 
$i$th component of $\tilde\nabla \vg$ at $x$ is equal to $
\left|\matrix g_1(x) & g_2(x)\\  \partial_ig_1(x) &
\partial_ig_2(x)\endmatrix \right|$, 
that is, to the signed area of the parallelogram spanned by 
$\vg(x)$ and $\partial_i\vg(x)$. Another interpretation: if one
represents $\vg(x)$ in polar coordinates, i.e.~via functions $\rho(x)$
and $\theta(x)$, it is straightforward to verify that 
$\tilde\nabla \vg(x)$ can be written as $\rho^2(x) \nabla \theta(x)$.

Loosely speaking, the
skew-gradient measures how different are the two functions from being
proportional to each other: it is easy to see that $\tilde\nabla \vg$
is identically equal to zero on an open set iff $g_1$ and $g_2$ are
proportional (with a locally constant coefficient). Therefore if the
image $\vg(V)\subset \br^2$ does not 
look like a part of a straight line passing through the origin, one
should expect the values of $\tilde\nabla \vg$ to be not very
small. Moreover, if the map $\vg$ is polynomial of degree $\le k$, then
$\tilde\nabla \vg$ is a polynomial map of degree $\le 2k-2$; in
particular, its norm is \cag\ for some $C,\alpha$. The results of the
previous section suggest that the latter property should be shared by
maps which are ``close to polynomial'' in the sense of Lemma 3.4 (that
is, for families of functions with some uniformly non-vanishing partial
derivatives).

The goal of the section is to prove the following result:

\proclaim{Proposition 4.1}   Let $U$ be an open subset of $\br^d$, $x_0\in
U$, and let 
$\Cal G\subset C^l(U)$ be a  family of maps
$\vg:U\mapsto \br^2$ such that 
$$
\text{ the family }\{\nabla g_i\mid \vg = (g_1,g_2)\in \Cal G,\ i = 1,2\}\text{ is compact in }C^{l-1}(U)\,.\tag 4.1
$$
 Assume also that 
$$
\inf \Sb{\vg\in \Cal G }\\ {\vv\in\br^2_1} \endSb 
\sup_{|\beta| \le l}
|\vv\cdot\partial_\beta \vg(x_0)| > 0\tag 4.2
$$
(in other words, the partial derivatives of $\vg$ at $x_0$ of order up to
$l$ uniformly in $\vg\in \Cal G$ generate $\br^2$).
Then there exists  a
neighborhood $V\subset U$ of $x_0$
such that 

{\rm (a)} $\|\tilde\nabla \vg\|$ is 
$(2C_{d,{l}},\frac1{d(2l-1)})$-good on $V$ for every $\vg\in \Cal G$ (here
$C_{d,{l}}$ is as in Lemma 3.2);

{\rm (b)} for every 
neighborhood $B\subset V$ of $x_0$ there exists $\rho = \rho(\Cal G, B)$ 
such that 
$$
\sup_{x\in B}\|\tilde\nabla \vg(x)\| \ge \rho\quad\text{ for every
}\vg\in \Cal G\,.
$$
\endproclaim

To prove this proposition, we will use two lemmas below. Note that 
in this section for convenience we will switch to  the Euclidean norm
$\|\vx\| = \|\vx\|_e$.  

\proclaim{Lemma 4.2} Let $B\subset \br^d$ be a ball of radius $r$ and
let $\vg$ be a $C^1$ map $B\mapsto \br^2$. Take $x_0\in
B$ such that $a = \vg(x_0)\ne 0$, denote  the line 
connecting $\vg(x_0)$ and the origin by $\Cal L$,   and let $\delta =
\sup_{x\in B}\|\vg(x) - \vg(x_0)\|$ and $w = \sup_{x\in 
B}\dist\big(\vg(x),\Cal L\big)$. Then 
$$
\sup_{x\in B}\|\tilde\nabla \vg(x)\| \ge
\frac{w(a-\delta)^2}{2r\sqrt{w^2 + (a+\delta)^2}}\,.
$$  
\endproclaim

\demo{Proof} Let us use polar coordinates, choosing $\Cal L$ to be the
polar axis.  Take
$x_1\in \overline B$ such that $\dist\big(\vg(x),\Cal L\big) = w$; then
one has
$$
\theta(x_1) \ge \sin \theta(x_1)  = \frac{\dist\big(\vg(x),\Cal
L\big)}{\|\vg(x_1)\|}\ge  \frac{w}{\sqrt{w^2 + (a+\delta)^2}}\,.
$$
Denote by $J$ the straight line segment $[x_0,x_1]\subset B$, and by
$u$ the unit vector proportional to $x_1 - x_0$. Restricting $\vg$
to $J$ and using Lagrange's
Theorem, one can find $y$ between $x_0$ and $x_1$ such that $\theta(x_1)
= 
\frac{\partial\theta}{\partial u}(y)|J|$. Then one has
$
|u\cdot \tilde\nabla \vg(y)| =
\rho^2(y)|\frac{\partial\theta}{\partial u}(y)| \ge
(a-\delta)^2\frac{\theta(x_1)}{|J|}$
which completes the proof.
\qed\enddemo

\proclaim{Lemma 4.3} Let $B\subset \br^d$ be a ball of radius $1$, and
let $\vp = (p_1,p_2):B\mapsto \br^2$  be a polynomial map  
of degree $\le l$ such that
$$
\sup_{x,y\in B}\|\vp(x) - \vp(y)\| \le 2\tag 4.3a
$$
(the diameter of the image of $\vp$ is bounded from above), and
$$
\sup_{x\in B}\text{\rm dist}\big(\Cal L,\vp(x)\big) \ge 1/8 \text{ for any straight line }\Cal L\subset \br^2\tag 4.3b
$$
(that is, the ``width'' of
 $\vp(B)$ in any direction is bounded from below). 
Then:

{\rm (a)} there exists a constant  
$0 < \gamma < 1$ (dependent only on $d$ and $l$) such that
$$
\sup_{x\in B}\|\tilde\nabla \vp(x)\| \ge \gamma\big(1 + 
\sup_{x\in B}\|\vp(x)\|\big)\,;\tag 4.4a
$$

{\rm (b)} there exists  $M \ge 1$ (dependent only on $d$ and $l$) such that
$$
\sup_{x\in B,\,i = 1,2}\|\nabla p_i(x)\| \le M\,.\tag 4.4b
$$ 
\endproclaim  
 

\demo{Proof} Let $\Cal P$ be the set of polynomial maps 
$\vp:B\mapsto \br^2$ 
of degree $\le l$ satisfying (4.3ab) and such that  $\sup_{x\in B}\|\vp(x)\| \le 6$. We first prove that there  exists $\gamma > 0$ such that 
(4.4a) holds for any $\vp\in\Cal P$. Indeed, otherwise from the compactness of $\Cal P$ it follows that there exists $\vp\in\Cal P$ such that $\tilde\nabla \vp(x)$ is identically equal to zero. Clearly this can only happen when all coefficients of $\vp$ are proportional to each other, which contradicts (4.3b). 

Now assume that $\vp$ satisfies (4.3ab) and $a = \|\vp(y)\| >
6$ for some $y\in B$. Then one can apply the previous lemma to the map
$\vp:B\mapsto \br^2$  to get $\sup_{x\in B}\|\tilde\nabla \vp(x)\| \ge \frac{\frac18(a -
2)^2}{2\sqrt{\frac1{64} + (a +
2)^2}} \ge \frac{\frac1{16}(a -
2)^2}{\sqrt{4(a -
2)^2}} \ge \frac1{32}(a -
2) \ge \frac1{64}(a +
1)$,
which finishes the proof of part (a). It remains to observe that part (b) 
trivially follows from the compactness of the 
set of polynomials of the form $\nabla \vp(x)$ where $\vp(x)$
satisfies (4.3a) 
 and has degree $\le l$. 
\qed\enddemo

\demo{Proof of Proposition 4.1} 
Choose $0 < \delta < 1/8$ such that 
$$
2C_{d,{l}}N_d\delta^{\frac1{d(2l-1)(2l-2)}} \le 1\,.
$$ 
From (4.1) and (4.2)  it follows that there
exists a neighborhood $V$ of $x_0$ and a positive $c$ such that for
every $\vg\in \Cal G$ 
one has 
$$
\forall\,{\vv\in\br^2_1} \quad\exists\,u\in\br^d_1\text{ and }k \le
l\text{ such that }\inf_{x\in V}\big|\vv\cdot\frac{\partial^k \vg}{\partial
u^k}(x)\big| \ge c\,, \tag 4.5a 
$$
and
$$
\sup_{x,y\in V}\|\partial_\beta\vg(x) - \partial_\beta\vg(y)\|  \le \frac{\delta c \gamma}{8 M l^l (l+1)!}\text{
for all  multiindices }\beta
\text{ with }|\beta| = {l}\,.\tag 4.5b
$$

In view of
Lemma 3.7, to show (a) it suffices to prove the following:  given any ball
$B = B(x_0,r)\subset V$ and a  $C^l$ map 
$\vg:B\to \br^2$ such that inequalities (4.5ab) hold  
for all   $x,y\in B$, 
one has 
$$
\big|\{x\in B\bigm| \|\tilde\nabla \vg(x)\| < \vre\cdot{\sup_{x\in B}\|\tilde\nabla \vg(x)\|}\}\big| \le
2C_{d,{l}}\vre^{\frac1{d(2l-2)}} 
|B|
\quad\text{whenever }\vre \ge \delta\,.
\tag 4.6
$$

We will do this in several steps.

{\bf Step 1.} Note that conditions (4.5ab), as well as the function $\tilde\nabla \vg$, will not change if one replaces $\vg$ by $L\vg$ where $L$ is any rotation of the plane $(g_1,g_2)$. Thus one can choose the ``$g_1$-axis'' in such a way that it is parallel to the line connecting two most distant points of $\vg(B)$.
By  (4.5a)  there exist $1 \le k_1,k_2 \le n$ and
$u_1,u_2\in\br^d_1$ such that
$|\frac{\partial^k g_i}{\partial
u_i^k}(x)| \ge c$ for $i = 1,2$ and all $x\in B$.
If $s_i$ stands for  $\sup_{x,y \in
B}|g_i(x) - g_i(y)|$, $i = 1,2$, 
then it follows from Lemma 3.6 that 
$$
s_i\ge \frac{c}{{k_i}^{k_i}({k_i}+1)!}(2r)^{k_i} \ge
\frac{c}{l^l(l+1)!}(2r)^l\,.\tag 4.7 
$$

{\bf Step 2.}  Here we replace the functions $g_i(x)$ by
$\frac1{s_i}g_i\big(x_0 + rx\big)$, and the ball $B$ by the unit ball $B(0,1)$.
This way the function $\tilde\nabla \vg$ will be multiplied by a
constant, and the statement (4.6) that we need to prove will be left
unchanged. However the partial derivatives of order $l$ of the functions $g_i(x)$ will be multiplied by factors $\frac{r^l}{s_i}$. In view of (4.7), the inequality (4.5b) will then imply
$$
\sup_{x,y\in B}\|\partial_\beta\vg(x) - \partial_\beta\vg(y)\|
\le {\delta \gamma}/{8M}\tag 4.8 
$$
for all  multiindices $\beta$
 with $|\beta| = {l}$
(here and until the end of the proof, $B$ stands for $B(0,1)$). Note
also that it follows from the construction that $\vg(B)$ is contained
in a translate of the square $[-\frac12,\frac12]^2$, and that  
$
\sup_{x\in B}\text{\rm dist}\big(\Cal L,\vg(x)\big) \ge 1/2\sqrt{2}
$
for any straight line $\Cal L\subset \br^2$. 

{\bf Step 3.}  Here we introduce the $l$-th degree Taylor polynomial
$\vp(x)$ of $\vg(x)$ at $0$. Using (4.8) one can show that $\vp$
is $\frac{\delta \gamma}{8M}$-close to $\vg$ in the $C^1$ topology,
that is, 
$$
\sup_{x\in B}\|\vg(x) - \vp(x)\| \le {\delta
\gamma}/{8M}\quad\text{and}\quad 
\sup_{x\in B}\|\nabla g_i(x) - \nabla p_i(x)\| \le {\delta
\gamma}/{8M},\ i = 1,2\,. 
$$
It follows that 
conditions (4.3ab) are satisfied by $\vp$, and therefore, by Lemma
4.3, the inequalities (4.4ab) hold.

{\bf Step 4.}  Now let us compare the functions $\tilde\nabla \vg$ and $\tilde\nabla \vp$: one has
$$
\tilde\nabla \vg - \tilde\nabla \vp =  \big(g_1(\nabla g_2- 
\nabla p_2) - (g_2 - p_2)\nabla g_1\big) -  \big((g_1 - p_1) \nabla p_2 - p_2
(\nabla g_1- \nabla p_1)\big)
$$
therefore 
$$
\split
\|\tilde\nabla \vg(x) - \tilde\nabla \vp(x)\| &\le \frac{\delta
\gamma}{8M}\big(\sup_{x\in B}|g_1(x)| + \sup_{x\in B}\|\nabla g_1(x)\| +  
\sup_{x\in B}|p_2(x)| + \sup_{x\in B}\|\nabla p_2(x)\| \big) \\
&\le \frac{\delta \gamma}{4M}(\big(\sup_{x\in B}\|\vp(x)\| +
\sup_{x\in B}\|\nabla p_2(x)\| \big) + \frac{\delta \gamma}{8M})
\\ &\un{(4.4b)}\le \frac{3\delta\gamma}8\big(1 +  
\sup_{x\in B}\|\vp(x)\|\big) \un{(4.4a)}\le
\frac{3\delta}{8}\sup_{x\in B}\|\tilde\nabla \vp(x)\|\,. 
\endsplit
$$

{\bf Step 5.} Finally we are ready to prove (4.6):  take $\vre$
between $\delta$ and $1$, 
put $s \df \sup_{x\in B}\|\tilde\nabla \vp(x)\|$ and observe that,
in view of Step 4, the set in the left hand side of (4.6) is contained
in  
$$
\left\{x\in B\bigm| \|\tilde\nabla \vp(x)\| - 
\tfrac{3\delta}{8}s< \vre(1 + \tfrac{3\delta}{8})s\right\}
 = \left\{x\in B\bigm| \|\tilde\nabla \vp(x)\| < \left(\vre +
\tfrac{3\delta}{8} 
(1 + \vre)\right)s\right\}\,.
$$
 Since $\vre + \tfrac{3\delta}{8} 
(1 + \vre) \le \vre + \tfrac{3\delta}{4} \le 2\vre$, and since
$\tilde\nabla \vp$ is a polynomial of degree not greater than $2l -
2$, one can apply Lemma 3.2 and 
conclude that the left hand side of (4.6) is not greater than 
$$
\big|\{x\in B\bigm| \|\tilde\nabla \vp(x)\| < 2\vre s\}\big| \le
C_{d,l} (2\vre)^{\frac1{d(2l-2)}}|B| \le 2C_{d,l}\vre^{\frac1{d(2l-2)}}|B|\,, 
$$
which 
finishes the proof of part (a). 

As for part (b), take a ball 
$B\subset V$ of radius $r$, and denote by $\hat B$ the ball with the
same center and 
twice smaller radius. It is
clear that there exists $\tau > 0$ such that for any $\vg\in \Cal G$
one can choose $y\in \hat B$ with  
$\|\vg(y)\| \ge \tau$ (otherwise, by a
compactness argument similar to that of Lemma 4.3, one would get that
$0|_{\hat B}\in \Cal G$, contradicting (4.2)). Also take $K \ge \tau/r$ such that  
$$
\sup_{\vg\in \Cal G,\,x\in B,\,u\in\br^d_1}\|\frac{\partial\vg}{\partial
u}(x)\| \le K\,.\tag 4.9 
$$

Now  let $B'\subset B$ be a ball
of radius $\tau/2K \le r/2$ 
centered at $y$. Take  $\vv\in\br^2_1$ orthogonal to
$\vg(y)$.  Applying Lemma 3.6  to $B'$ and the function $\vv\cdot\vg(x)$
one gets $\sup_{x\in
B'}|\vv\cdot\vg(x)|\ge \frac{c}{k^k(k+1)!}(\tau/K)^k \ge
\frac{c}{l^l(l+1)!}(\tau/K)^l$. On the other hand  (4.9) shows that  $\sup_{x\in
B'}\|\vg(x) - \vg(y)\|$ is not greater than $\tau/2$. Now one can apply
Lemma 4.2 to the map
$\vg:B'\mapsto \br^2$  to get 
$$
\sup_{x\in B'}\|\tilde\nabla \vg(x)\| \ge
\frac{c}{l^l(l+1)!}\Big(\frac\tau K\Big)^{l-1}
\frac{(\tau/2)^2}{\sqrt{\big(\frac{c}{l^l(l+1)!}(\tau/K)^l\big)^2
+ (3\tau/2)^2}}\,,
$$
giving a uniform lower bound for $\sup_{x\in B}\|\tilde\nabla
\vg(x)\|$.
\qed\enddemo

\heading{5. Theorem 1.4 and lattices}
\endheading 

Roughly speaking, the method of lattices simply allows one to write
down the system of inequalities (1.7b) from Theorem 1.4 in an
intelligent way. In what follows, we let $m$ stand for
${n+d+1}$. Denote  the   
standard basis of  
$\br^{m}$ by $\{\ve_0,\ve_1^*,\dots,\ve_d^*,\ve_1,\dots,\ve_n\}$.
Also denote by $\Lambda$ the intersection of  
$\bz^{m}$ with the span of $\ve_0,\ve_1,\dots,\ve_n$,  that is,
$$
\Lambda = \left\{\left.\left(\matrix
p  \\ 
0 \\ \vq
\endmatrix \right)\right| p\in\bz,\ \vq\in\bz^n\right\}\,.\tag 5.1
$$ 
Take $\vf:U\mapsto \br^n$ is as in Theorem 1.4 and let $U_x$ stand for
the matrix  
$$
U_x \df \left(\matrix
1 & 0 & \vf(x)  \\
0 & I_d & \nabla\vf(x) \\0 & 0 & I_n
\endmatrix \right)\in SL_{m}(\br)\,.\tag 5.2
$$
Note that 
$
U_x\left(\matrix
p  \\ 
0 \\ \vq
\endmatrix \right) = \left(\matrix \vf(x)\vq + p \\ \nabla\vf(x)\vq
\\ \vq \endmatrix \right)$
is the vector whose components appear in the right hand sides of
the inequalities (1.7b). Therefore the fact that 
 there exists $\vq\in\bz^n\nz$ satisfying (1.7b)
implies  the existence of a nonzero element  of $U_x\Lambda$
which belongs to a certain parallelepiped in $\br^{m}$. Our strategy
will be as follows: 
we will find a diagonal matrix $D \in GL_{m}(\br)$ which
transforms the above 
parallelepiped into a small cube; then the solvability of the above
system of inequalities will force
the lattice 
$DU_x\Lambda$ to have a small nonzero vector, and we will  use
a theorem proved by methods from \cite{KM1} (see Theorem 6.2 below)
to estimate the measure of the set 
of $x\in B$  for which it can happen. 

Specifically, take $\delta, K, {T}_1,\dots,{T}_n$ as in Theorem 1.4,
fix $\vre > 0$ and    denote 
$$
D =
\text{diag}(a_0^{-1},a_*^{-1},\dots,a_*^{-1},a_1^{-1},\dots,a_n^{-1})\,,\tag
5.3 
$$ 
where
$$
a_0 = \frac{\delta}\vre,\ a_* =  \frac{K}\vre,\ a_i =
\frac{{T}_i}\vre, \ i = 
1,\dots,n\,.\tag 5.4
$$
It can be easily seen that  the set (1.7b) is exactly equal to
$$
\big\{x\in B\mid \|DU_x\vv\| < \vre \text{ for some
}\vv\in\Lambda\nz \big\}\,,\tag 5.5
$$
where $\|\cdot\|$ stands for the supremum norm. However from this
point on it will be more convenient to use the Euclidean norm
$\|\cdot\|_e$ on
$\br^{m}$. Let us now state a theorem from which Theorem 1.4 can be
easily derived.  

\proclaim{Theorem 5.1} Let $U$, $x_0$, $d$, $l$, $n$ and
$\vf$
be as in Theorem 1.4. Take  $\Lambda$ as in {\rm (5.1)} and
$U_x$ as in {\rm (5.2)}. Then there exists a neighborhood $V\subset U$ of
$x_0$  with the following property: for
any ball $B\subset V$ there exists
$E > 0$ such that for any diagonal matrix $D$ as in {\rm (5.3)} with 
$$
0 < a_0 \le 1,\quad a_n \ge\dots \ge a_1\ge 1\quad  \text{and}\quad 0
<  a_* \le (a_0a_1\dots  a_{n-1})^{-1}
\tag 5.6
$$ 
and any positive $\vre$, one has 
$$
\left|\big\{x\in B\mid \|DU_x\vv\|_e < \vre \text{ for some
}\vv\in\Lambda\nz \big\}\right| \le E
\vre^{\frac1{d(2l-1)}} |B| 
\,.\tag 5.7
$$
\endproclaim

\demo{Proof of Theorem 1.4 modulo Theorem 5.1} Take $V\subset U$ and,
given any ball $B\subset V$, choose $E$ as in the above theorem. Then
take $\delta, {T}_1, \dots, {T}_{n}$ and
$K $ satisfying (1.7a). Observe that without
loss of generality one can assume that ${T}_1\le  \dots \le
{T}_{n}$. (Otherwise one can replace ${T}_1, \dots, {T}_{n}$  by a
permutation ${T}_{i_1},
\dots, {T}_{i_n}$ with  ${T}_{i_1}\le 
\dots \le {T}_{i_n}$, and consider the $n$-tuple
$(f_{i_1},\dots,f_{i_n})$ instead of the original one, which will
still be nondegenerate at $x_0$.)

Now  take $\vre$ as in (1.7c) 
and define
$a_0,a_*,a_1,\dots,a_n$ as in 
(5.4).  Then all the constraints (5.6) easily follow (indeed, (1.7c)
shows that $a_0 \le 1$ and $a_0a_*a_1\dots  a_{n-1} \le 1$, while
(1.7a) implies that $\vre \le 1$, hence $a_i \ge 1$).  
Thus  Theorem 5.1  applies, and to complete the proof it remains to
observe that the set (1.7b) $=$ (5.5) is contained in 
$$
\big\{x\in B\bigm| \|DU_x\vv\|_e < \vre\sqrt{m} \text{ for some
}\vv\in\Lambda\nz \big\}\,,
$$
hence its measure is not greater than 
$E m^{\frac1{2d(2l-1)}}
\vre^{\frac1{d(2l-1)}} |B|$. \qed \enddemo

\heading{6. Lattices and posets}
\endheading 

The proof of Theorem 5.1 will depend on a result from \cite{KM1}
involving mappings of 
partially ordered sets into spaces of \cag\ functions. Let us recall
some terminology from \cite{KM1, \S 4}. 
For $d\in\bn$, $k\in \bz_{\sssize +}$ and $C,\alpha,\rho   > 0$, define $\ca$ to
be the set of triples $\sfb$ where $S$ is a partially ordered set
({\sl poset\/}),  $B = B(\vx_0,r_0)$, where $\vx_0\in\br^d$ and $r_0 >
0$,  and $\ph$ is a mapping from $S$ to the space of 
continuous functions on $\tilde B \df B\big(\vx_0,3^kr_0\big)$ (this
mapping will be denoted by $s\to 
\ph_s$) such that the following holds:


\roster
\item"(A0)" the length of $S$ is not greater than $k$;

\item"(A1)" $\forall\,s\in S\,,\quad \ph_s$ is \cag\ on $B$;

\item"(A2)" $\forall\,s\in S\,,\quad\|\ph_s\|_{\sssize B} \ge \rho $.

\item"(A3)"  $\forall\,x\in B,\quad\#\{s\in S\bigm| |\ph_s(x)| <
\rho\} < \infty$. 
\endroster

Then, given $\sfb\in\ca$ and $\vre > 0$, say that a point $x\in B$ is {\sl $(\vre,S,\ph)$-marked\/} if there exists a linearly ordered subset $\Sigma_x$ of $S$ such that

\roster
\item"(M1)" $\vre\le |\ph_s(x)| \le \rho \quad \forall\,s\in\Sigma_x$;

\item"(M2)" $|\ph_s(x)| \ge \rho \quad\forall\,s\in S \ssm \Sigma_x$
comparable with any element of $\Sigma_x$. 
\endroster

We will denote by $\p(\vre,S,\ph,B)$ the set of all the $(\vre,S,\ph)$-marked points $x\in B$.

 \proclaim{Theorem 6.1 \rm (cf.~\cite{KM1, Theorem 4.1})} Let
$d\in\bn$, $k\in \bz_{\sssize +}$ and $C,\alpha,\rho  > 0$ be
given. Then  for all $\sfb\in\ca$ and $ 0 < \vre \le \rho$ one has
$$
\left|B\smallsetminus \p(\vre,S,\ph,B)\right| \le kC\big(3^dN_d\big)^k  \left(\frac\vre \rho \right)^\alpha  |B|\,.
$$
\endproclaim

We will apply
Theorem 6.1 to the poset of subgroups of the group of integer points
of a finite-dimensional real vector space $W$. 
For a  discrete subgroup $\Gamma$  of 
$W$, we will denote by 
$\Gamma_\br$ the minimal linear subspace of $W$ containing
$\Gamma$. Let  $k =  \text{dim}(\Gamma_\br)$ be the {\sl rank\/} of
$\Gamma$; say that $\vw\in 
\bigwedge^k(W)$ {\sl represents\/} $\Gamma$ if  
$$
\vw =  \cases &1\qquad\qquad\quad\text{\ \ \ if } k = 0\\
&\vv_{1}\wedge\dots\wedge \vv_{k}\quad\text{if }k > 0 \text{ and
}\vv_{1},\dots, \vv_{k} \text{ is a basis of }\Gamma\,.\endcases 
$$
We will need this  exterior power representation mainly to be able to
measure the ``size'' of discrete subgroups. Namely, these ``sizes''
will be given by suitable ``norm-like'' functions $\nu$ on the exterior algebra 
$\bigwedge(W)$ of $W$, and we will set 
$$
\nu(\Gamma) = \nu(\vw) \quad\text{if }\vw\text{ represents }\Gamma\,.\tag 6.1
$$
More precisely, let us say that a function $\nu:
\bigwedge(W)\mapsto \brp$ is {\sl submultiplicative\/} if

\roster
\item"(i)"  $\nu$ is continuous (in the natural topology);

\item"(ii)" it is homogeneous, that is, $\nu(t \vw) =
|t|\nu(\vw)$ for all $t\in\br$ and $\vw\in
\bigwedge(W)$; 

\item"(iii)" $\nu(\vu\wedge \vw) \le \nu(\vu)\nu(\vw)$ for all $\vu,\vw\in
\bigwedge(W)$.
\endroster

Note that in view of (ii), (6.1) is a correct definition of
$\nu(\Gamma)$. 

Examples: if $W$ is a Euclidean space, one can  extend the Euclidean
structure to $\bigwedge(W)$ (by making $\bigwedge^i(W)$ and
$\bigwedge^j(W)$ orthogonal for $i\ne j$); clearly then the  Euclidean
norm $\nu(\vw) = 
\|\vw\|$ is submultiplicative. In this case the restriction of $\nu$ to
$W$ coincides with the usual (Euclidean) norm  on $W$. Also if $\Cal
W\subset \bigwedge(W)$ is 
an ideal, one can define $\nu(\vw)$ to be the norm of the projection
of $\vw$ orthogonal to $\Cal W$. If this ideal is orthogonal to
$W\subset \bigwedge(W)$, again the function $\nu$ will  coincide with
the norm when restricted to $W$. 

We will also need the notion of {\sl primitivity\/} of a discrete
subgroup. If $\Lambda$ is a discrete subgroup  of
$W$,   say that a
subgroup $\Gamma$ of $\Lambda$ is {\sl primitive\/}  (in $\Lambda$) if
$\Gamma = \Gamma_\br\cap \Lambda$, and denote by $\Cal L(\Lambda)$ the
set of all nonzero  primitive subgroups of $\Lambda$.    Example: a
cyclic subgroup of 
$\bz^l$   is primitive in $\bz^l$  iff it is generated by a {\sl 
primitive\/} vector (that is, a vector which is not equal to a
nontrivial multiple of another element of $\bz^l$). Note that the
inclusion relation makes $\Cal L(\Lambda)$ a poset, its length being
equal to the rank of $\Lambda$. 

The following result, which we will derive from Theorem 6.1,   is a
generalization of 
Theorem 5.2 from \cite{KM1}.

 \proclaim{Theorem 6.2} Let $W$ be a finite-dimensional real vector
space, $\Lambda$ a discrete subgroup of $W$ of rank $k$,  and let a
ball $B = B(x_0,r_0)\subset \br^d$ and a map $H:\tilde B \to GL(W)$
be given, where $\tilde B$ stands for $B(x_0,3^kr_0)$. Take
$C,\alpha > 
0$, $0 < \rho  
\le 1$, and let $\nu$ be a submultiplicative function on $\bigwedge(W)$.
Assume that for any $\Gamma\in\Cal L(\Lambda)$,

{\rm(i)} the function $x\mapsto \nu\big(H(x)\Gamma\big)$ is \cag\ on
$\tilde B$,  and  

{\rm(ii)} $\exists\,x\in B$ such that 
$\nu\big(H(x)\Gamma\big) \ge \rho$. 

\noindent Also assume that

{\rm(iii)}  $\forall\,x\in \tilde B,\quad\#\big\{\Gamma\in\Cal L(\Lambda)\bigm| \nu\big(H(x)\Gamma\big) <
\rho\big\} < \infty$.

\noindent Then 
 for any  positive $ \vre \le \rho$ one has 
$$
\left|\big\{x\in B\bigm| \nu \big(H(x)\vv\big) < \vre \text{ for some
}\vv\in\Lambda\nz \big\}\right| \le k(3^dN_d)^k\cdot C \left(\frac\vre
\rho \right)^\alpha  |B|\,.\tag 6.2 
$$
\endproclaim
 
Note that when $\nu|_W$ agrees with the Euclidean norm, (6.2) estimates the
measure of $x\in B$  
for which the subgroup $H(x)\Lambda$ has a nonzero vector with
length less than $\vre$. 

\demo{Proof} We will apply Theorem 6.1 to the triple
$\sfb$, where  $S = \Cal 
L(\Lambda)$ and  $\ph$ is
defined by $\ph_{\sssize \Gamma}(x) \df \nu\big(H(x)\Gamma\big)$. It is easy to
verify that  $\sfb\in\ca$.  
Indeed, the functions $\ph_{\sssize \Gamma}$ are continuous since so
is $H$ and $\nu$,  property (A0) is clear, (A1) is given by (i), (A2) by
(ii) and (A3) by (iii).   

In view of Theorem 6.1, it remains to prove that 
$$
\p(\vre,S,\ph,B) \subset  \big\{x\in B\mid \nu\big(H(x)\vv\big) \ge \vre\text{ for all }\vv\in\Lambda\nz\big\}\,.\tag 6.3
$$
Take an $(\vre,S,\ph)$-marked point $x\in B$, and let $\{0\} =
\Gamma_0 \subsetneq \Gamma_1 \subsetneq\dots \subsetneq\Gamma_m =
\Lambda$ be all the elements of $\Sigma_x \cup
\big\{\{0\},\Lambda\big\}$. Take any $\vv\in\Lambda\nz$. Then there
exists $i$, $1 \le i \le m$, such that $\vv\in \Gamma_i\ssm
\Gamma_{i-1}$. Denote $(\Gamma_{i-1} + \br\vv)\cap \Lambda$ by
$\Delta$. Clearly $\Delta$ is a primitive subgroup of $\Lambda$
contained in $\Gamma_i$, therefore comparable to any element of
$\Sigma_x$. By submultiplicativity of $\nu$ one has
$\nu\big(H(x)\Delta\big) \le
\nu\big(H(x)\Gamma_{i-1}\big)\nu\big(H(x)\vv\big)$. Now one can use
properties (M1) and 
(M2) to deduce 
that 
$
|\ph_{\sssize \Delta}(x)| = \nu\big(H(x)\Delta\big) \ge \min(\vre, \rho) = \vre\,,
$
and then conclude that 
$$
\nu\big(H(x)\vv\big) \ge
{\nu\big(H(x)\Delta\big)}/{\nu\big(H(x)\Gamma_{i-1}\big)}  \ge \vre/\rho  \ge \vre\,.
$$
This shows (6.3) and completes the proof of the theorem. \qed\enddemo

\heading{7. Proof of Theorem 5.1} \endheading

Here we take $U$, $x_0$, $d$, $l$, $n$ and
$\vf$ as in Theorem 1.4, set 
 $m = n+d+1$ and $W = \br^{m}$,
 and  use the notation introduced in \S 5. 
Denote by $W^*$ the $d$-dimensional subspace spanned by
$\ve_1^*,\dots,\ve_d^*$, so that $\Lambda$ as in (5.1) is equal to the
intersection of  
$\bz^{m}$ and $(W^*)^\perp$. Also let $\Cal H$ be the family of functions 
$H:U\mapsto GL_m(\br)$ given by 
$H(x) = DU_x
$,
where  $U_x$ is as in (5.2), $D$ as in (5.3) with coefficients 
satisfying (5.6).
 
In order to use Theorem 6.2, we also need to choose the submultiplicative
function $\nu$ on $W$ 
in a special way.  Namely, we let $\Cal W\subset \bigwedge(W)$ be the
ideal generated by  $\bigwedge^2(W^*)$,  denote by $\pi$ the
orthogonal projection with kernel 
$\Cal W$, and take $\nu(\vw)$ to be the
Euclidean norm of $\pi(\vw)$. In other words, if $\vw$ is written as a
sum of exterior products of base vectors $\ve_i$ and $\ve^*_i$, to
compute $\nu(\vw)$ one
should ignore components containing $ \ve^*_i \wedge \ve^*_j$, $1\le i
\ne j \le d$, and take the norm of the sum of the remaining
components.

Since $\nu|_W$ agrees with the
Euclidean norm, to derive Theorem 5.1 from Theorem 6.2 it suffices to
find a neighborhood $\tilde V\ni x_0$ such that  

\roster 
\item"{$\boxed{1}$}" there exists $C>0$ such that all the functions 
$x\mapsto \nu\big(H(x)\Gamma\big)$, where  $H\in \Cal H$ and $\Gamma\in\Cal L(\Lambda)$, are $(C,\frac1{d(2l-1)})$-good on
$\tilde V$,

\item"{$\boxed{2}$}" for every ball $B\subset \tilde V$ there exists $\rho > 0$  such that 
$\sup_{x\in B}\nu\big(H(x)\Gamma\big) \ge \rho$ for all $H\in \Cal H$ and $\Gamma\in\Cal L(\Lambda)$, and

\item"{$\boxed{3}$}"  for all $x\in \tilde V$ and $H\in \Cal H$ one has $\#\big\{\Gamma\in\Cal L(\Lambda)\bigm| \nu\big(H(x)\Gamma\big) \le
1\big\} < \infty$.
\endroster

Indeed, then one can take a smaller neighborhood $V$ of $x_0$ such that whenever $B = B(x,r)$ lies in $V$, its dilate $\tilde B = B(x,3^{n+1}r)$ is contained in $\tilde V$. This way it would follow from Theorem 6.2 that for any $B\subset V$ the measure of the
set in (5.7) is not greater than $C (n+1)(3^dN_d)^{n+1}\left(\vre
/\rho \right)^{1/d(2l-1)}  |B|$ for any $\vre \le \rho$, therefore 
not greater than 
$$
\max\big(C(n+1)(3^dN_d)^{n+1},1\big)\rho^{-1/d(2l-1)} \vre
^{1/d(2l-1)}  |B|
$$ for any positive $\vre$. 

Thus we are led to explicitly computing the functions
$\nu\big(H(x)\Gamma\big)$ for arbitrary choices of subgroups 
$\Gamma\subset \Lambda$ and positive numbers $a_i$, $i =
0,*,1,\dots,n$. In fact, we will be doing it in two different ways,
which will be relevant for checking conditions $\boxed{1}$ (along with $\boxed{3}\,$) and $\boxed{2}$
respectively. 
 
Let $k$ be the rank of $\Gamma$. The claims are trivial for $\Gamma =
\{0\}$, thus we can set $1\le k \le n+1$.   
Since  $D\Gamma_\br$ is a $k$-dimensional subspace of
$(W^*)^\perp = \br\ve_0\oplus\br\ve_1\oplus\dots\oplus\br\ve_n$, it 
 is possible to choose an orthonormal set $\vv_1,\dots,\vv_{k-1}\in
D\Gamma_\br$ such that  
each $\vv_i$, $i = 1,\dots,k-1$, is orthogonal to $\ve_0$. Now let us
consider two cases:

{\bf Case 1.}  $D\Gamma_\br$ contains $\ve_0$; then $\{\ve_0,\vv_1,\dots,\vv_{k-1}\}$ is a basis of $\br\ve_0\oplus D\Gamma_\br$. Thus one can find    $\vw\in\bigwedge^k(\br^{m})$   representing $\Gamma$ such that $D\vw$ can be written  
as $a \ve_{0}\wedge\vv_{1}\dots\wedge \vv_{k-1}$ for some $a > 0$.

{\bf Case 2.} $D\Gamma_\br$ does not contain $\ve_0$; then it is possible 
to choose $\vv_0\in \br\ve_0\oplus D\Gamma_\br$ such that $\{\ve_0,\vv_0, \vv_1,\dots,\vv_{k-1}\}$ is  an orthonormal  basis of $\br\ve_0\oplus D\Gamma_\br$. In this case,
one can represent $\Gamma$ by $\vw$ such that 
$$
D\vw = (a \ve_{0}  + b\vv_{0})\wedge\vv_{1}\dots\wedge \vv_{k-1} =  a \ve_{0}\wedge\vv_{1}\dots\wedge \vv_{k-1} + b
\vv_{0}\wedge\vv_{1}\wedge\dots\wedge \vv_{k-1}\tag 7.1 
$$
for some $a,b\in \br\nz$. In fact (7.1) is valid in
Case 1 as well; one simply has to put $b$ equal to zero and bear in mind
that the vector $\vv_0$ is not defined.

Now write
$$
H(x)\vw = DU_x D^{-1}(D\vw)\,,\tag 7.2
$$
and introduce the $m$-tuple (interpreted as a row vector) of functions 
$$
\hat \vf(x) =
\big(1,0,\dots,0,\frac{a_1}{a_0}f_1(x),\dots,\frac{a_n}{a_0}f_n(x)\big)\,.
$$
It will also be convenient to ``identify $\br^d$ with $W^*$'' and
introduce the $W^*$-valued gradient 
$
\nabla^* = \sum_{i = 1}^d \ve^*_i\partial_i
$ of a scalar function on $U\subset\br^d$, so that $\nabla^*f(x) = \sum_{i =
1}^d \partial_if(x)\ve^*_i$. As a one step further, we will let the
skew-gradients discussed in \S 4 take values in $W^*$ as well. That
is, for a map $\vg = (g_1,g_2):\br^d\mapsto \br^2$ we
will define
$\tilde\nabla^* \vg: \br^d\mapsto
W^*$ by
$$
\tilde\nabla^* \vg(x)\df g_1(x) \nabla^*g_2(x) - g_2(x)
\nabla^*g_1(x)\,.
$$

Then it becomes straightforward to verify that $DU_x D^{-1}\ve_0 = \ve_0$ and 
$$
DU_x D^{-1}\vv = \vv + \big(\hat \vf(x)\vv \big)\ve_0 +
{}\frac{a_0}{a_*} \nabla^*\big(\hat \vf(x)\vv\big)
$$
 whenever $\vv$ is orthogonal to $\ve_0$ and $W^*$.
Therefore the space $\br\ve_0\oplus W^*\oplus
D\Gamma_\br$ 
is invariant under $DU_x D^{-1}$; hence we can restrict
ourselves to the coordinates of $DU_x D^{-1}$-image of $
D\vw$ with respect to the basis chosen above, and write
$$
DU_x D^{-1}(\ve_{0}\wedge\vv_{1}\dots\wedge \vv_{k-1}) = \ve_{0}\wedge\vv_{1}\wedge\dots\wedge \vv_{k-1} + \dsize{}\frac{a_0}{a_*}\sum_{i = 1}^{k-1} \pm  \nabla^*\big(\hat \vf(x)\vv_i\big)\,
\ve_{0}\wedge 
\bigwedge_{s \ne i}\vv_s + \vw^*_1 
$$
and
$$
\aligned
DU_x D^{-1}(\vv_{0}\wedge\dots\wedge \vv_{k-1})\  = \ \vv_{0}\wedge\dots\wedge \vv_{k-1} \ +\  \dsize\sum_{i = 0}^{k-1}
\pm \big(\hat \vf(x)\vv_i \big)\, \ve_0 \wedge\bigwedge_{s \ne
i}   
\vv_s \\ +\  \dsize{}\frac{a_0}{a_*}\sum_{i = 0}^{k-1} \pm \nabla^*\big(\hat
\vf(x)\vv_i\big)\wedge 
\bigwedge_{ s \ne i} \vv_s \ 
+ \ \dsize{}\frac{a_0}{a_*}\sum\Sb{i = 0 }\\ {j > i}\endSb ^{k-1} \pm\tilde\nabla^*\left( 
\hat \vf(x)\vv_i,\hat \vf(x)\vv_j\right)\wedge\ve_{0}\wedge 
\bigwedge_{s \ne i,j}\vv_s + \vw^*_2 \,, \endaligned 
$$
where $\vw^*_1$ and $\vw^*_2$ belong to 
$\bigwedge^2(W^*)\wedge\bigwedge(W)$.

Collecting terms and using (7.1) and (7.2), one finds  that
$$
\aligned
\pi\big(H(x)\vw\big) \ &= \ \big(a + b\hat \vf(x)\vv_0\big)\, \ve_{0}\wedge\vv_{1}\wedge\dots\wedge \vv_{k-1} \ + \ b
\,\vv_{0}\wedge\dots\wedge \vv_{k-1}\\ 
&+\  b\dsize\sum_{i = 1}^{k-1}
\pm  \big(\hat\vf(x)\vv_i\big)\, \ve_0 \wedge\bigwedge_{s \ne
i}   
\vv_s \ + \ b\dsize{}\frac{a_0}{a_*}\sum_{i = 0}^{k-1} \pm \nabla^*\big(\hat
\vf(x)\vv_i\big)\, \wedge 
\bigwedge_{ s \ne i} \vv_s \\
&+\  \dsize{}\frac{a_0}{a_*}\sum_{i = 1}^{k-1} \pm 
\tilde\nabla^*\big(\hat
\vf(x)\vv_i, a + b\hat \vf(x)\vv_0 \big)\,\wedge 
\ve_{0}\wedge 
\bigwedge_{s \ne 0,i}\vv_s \\
&+ \ b\dsize{}\frac{a_0}{a_*}\sum\Sb{i,j = 1 }\\ {j > i}\endSb ^{k-1}
\pm \tilde\nabla^*\big(\hat \vf(x)\vv_i,
 \hat \vf(x)\vv_j\big)\,\wedge \ve_{0}\wedge 
\bigwedge_{s \ne i,j}\vv_s \endaligned \tag 7.3
$$

Now the stage is set to prove condition $\boxed{1}\,$. Indeed, in view of Lemma
3.1(bc), it suffices to show that the norms of each of the summands in
(7.3) are \cag\ functions. Note that the elements appearing 
in the first two lines are  linear
combinations of either $1,f_1,\dots,f_n$ or $\nabla^* f_1,\dots,\nabla^*
f_n$, hence one can, 
using  Corollary
3.5, find $\tilde V_1\ni x_0$ and $C_1>0$ such that all these norms are
$(C_1,1/n)$-good on $\tilde V_1$.  The same can be said about the rest of the
components when $b = 0$. Otherwise   one can observe 
that the summands in the  last two lines of (7.3)
are of the form $\pm \tilde\nabla^*(L\vg)$ where $L$ is some linear
transformation  
of $\br^2$ and $\vg$ belongs to 
$$
\Cal G \df \big\{\big(\vf(x)\vu_1, \vf(x)\vu_2   + u_0\big)\mid u_0\in\br, \vu_1 \perp \vu_2 \in \br^n_1 \big\}\,.\tag 7.4
$$
 Since $\tilde\nabla^*(L\vg)$ is proportional to $\tilde\nabla^*\vg$,
it suffices, 
in view of Lemma 3.1(a), to prove the existence of $C_2 \ge C_1$ and a
neighborhood  
$\tilde V_2 \subset \tilde V_1$ of $x_0$ such that the norms of
$\tilde\nabla^*\vg$  are 
$(C_2,\frac1{d(2l-1)})$-good on $\tilde V_2$ for any $\vg\in \Cal G$.
The latter is 
a direct consequence of Proposition 4.1(a), since the family $\{\nabla g_i\mid
\vg = (g_1,g_2)\in \Cal G\}$ is compact in $C^{l-1}(U)$,  
and it follows from the nondegeneracy of $\vf$ at $x_0$ that condition (4.2)
is satisfied. 

One can also use (7.3) to prove  $\boxed{3}\,$. Indeed, looking at the first line of (7.3) one sees that 
$$\nu\big(H(x)\Gamma\big)\le \max\big(|a + b\hat \vf(x)\vv_0|, |b|\big)\,,
$$ thus $\nu\big(H(x)\Gamma\big) \le 1$ would imply $|b|\le 1$ and $|a| \le 1 + \|\hat \vf(x)\|$.  On the other hand, in view of (7.1), the norm of $D\vw$ is equal to $\sqrt{a^2 + b^2}$,  and  from the discreteness of $\bigwedge(D\Lambda)$ in $\bigwedge(\br^m)$ it follows that for any $R > 0$ the set of $\Gamma\in\Cal L(\Lambda)$  such that both $|a|$ and $|b|$ are bounded from above by $R$ is finite.

We now turn to condition $\boxed{2}\,$. Take any neighborhood $B$ of $x_0$. It
follows from the 
linear independence of the functions $1,f_1,\dots,f_n$ and from the
linear independence of their gradients 
that there exists $\rho_1 > 0$ such that 
$$
\forall\,\vv\in\br^n_1\ \forall\,v_0\in\br\text{ one has }\sup_{x\in
B}|\vf(x)\vv 
 + v_0| \ge \rho_1 \text{ and } 
\sup_{x\in B}|\nabla\big(\vf(x)\vv\big)| \ge \rho_1 \,.\tag 7.5
$$
Also let  $\rho_2 = \rho(\Cal G,B)$ where $\Cal G$ is the class of $2$-tuples of functions defined in (7.4), and
$\rho(\Cal G,B)$ is as in Proposition 4.1(b). Consider
$$
M \df \max\big(\sup_{x\in B}\|\vf(x)\|,\sup_{x\in B}\|\nabla\vf(x)\|\big)\,;
$$
we will show that (ii) will hold for any nonzero subgroup $\Gamma$ of
$\Lambda$ if we choose  
$$
\rho  = \frac{\rho_1\rho_2}{\sqrt{\rho_1^2 + (\rho_2 + 2M^2)^2}}\,.
$$

First let us consider the case $k = \text{dim}(\Gamma_\br) = 1$. Then
$\Gamma$ can be represented by a vector $\vv =
(v_0,0,\dots,0,v_1,\dots,v_n)\t$ with integer coordinates, and it is
straightforward to verify that the first coordinate of 
$H(x)\vv$ will be equal to $\frac1{a_0}\big(v_0 + v_1f_1(x) + \dots +
v_nf_n(x)\big)$, which will deviate from zero by not less than $\rho$ at
some point of $B$ due to (7.5) and since $\rho \le \rho_1$ and $a_0 \le 1$. 

Now let $k$ be greater than $1$. Our method will be similar to that of
the proof of the  previous lemma: given $\vw\in\bigwedge^k(\br^{m})$ representing  $\Gamma$, we will choose a
suitable orthogonal decomposition  of $\bigwedge^k(\br^{m})$  and then show
that the norm of the projection of $H(x)\vw$ to some
subspace will be not less than $\rho$ for some $x\in B$. 

In order to prove the desired estimate, it is important to pay special 
attention to the vector $\ve_n$, which is the eigenvector of $D$ with the
smallest eigenvalue. We do it by first choosing an orthonormal set
$\vv_1,\dots,\vv_{k-2}\in \Gamma_\br$ such that  each $\vv_i$, $i =
1,\dots,k-2$, is orthogonal to both $\ve_0$ and $\ve_n$. Then choose
$\vv_{k-1}$ orthogonal to $\vv_i$, $i = 1,\dots,k-2$, and to $\ve_0$
(but in general not to $\ve_n$). Now, if necessary (see the remark
after (7.1)), 
choose a vector $\vv_0$ to complete $\{\ve_0,\vv_1,\dots,\vv_{k-1}\}$
to an orthonormal basis of $\br\ve_0\oplus \Gamma_\br$. This way,
similarly to (7.1), 
we will represent $\Gamma$ by $\vw$ of the form 
$$
\vw = (a \ve_{0}  + b\vv_{0})\wedge\vv_{1}\dots\wedge \vv_{k-1} =  a \ve_{0}\wedge\vv_{1}\dots\wedge \vv_{k-1} + b
\vv_{0}\wedge\dots\wedge \vv_{k-1}\tag 7.6
$$
for some $a,b\in \br$ with $a^2 + b^2 \ge 1$. As before, we will use (7.6) even when $\vv_0$ is not defined, in this case the coefficient $b$ will vanish.

Now, similarly to the proof of condition (i), introduce the  $m$-tuple
of functions  
$$
\check \vf(x) =
\big(1,0,\dots,0,f_1(x),\dots,f_n(x)\big)\,,
$$
 and observe that $U_x\ve_0 = \ve_0$ and 
$$
U_x\vv = \vv + \big(\check \vf(x)\vv \big)\ve_0 +
\nabla^*\big(\check \vf(x)\vv \big)
$$
 whenever $\vv$ is orthogonal to $\ve_0$ and $\ve_*$. Using this and
(7.6), one can obtain 
an expression analogous to (7.3). This time however we are interested
only in the terms of the form $\ve_{0}\wedge\ve^*_i\wedge \vw'$, where
$\vw'$ is orthogonal to $
\bigwedge^{k-2}\big((\br\ve_0 \oplus W^*)^\perp\big)$ (note that
these terms are present only if $k \ge 2$). Namely, let us write  
$$
\aligned
\pi(U_x\vw) \ &= \ \left(a + b\big(\check \vf(x)\vv_0\big)\right)\, \ve_{0}\wedge\vv_{1}\wedge\dots\wedge \vv_{k-1} \ + \ b
\,\vv_{0}\wedge\dots\wedge \vv_{k-1} \\
&+\  b\dsize\sum_{i = 1}^{k-1}
\pm \big(\check \vf(x)\vv_i\big)\, \ve_0 \wedge\bigwedge_{s \ne
i}   
\vv_s \ + \ b\sum_{i = 0}^{k-1} \pm \nabla^*\big(\check
\vf(x)\vv_i \big) \wedge\bigwedge_{ s \ne i} \vv_s  \ + 
\  \ve_{0}\wedge \check\vw(x)\,,
\endaligned 
$$
where
$$
\check\vw(x) \df \sum_{i = 1}^{k-1} \pm \tilde\nabla^*\big(\check
\vf(x)\vv_i, a + b\check \vf(x)\vv_0 \big)\wedge 
\bigwedge_{s \ne 0,i}\vv_s + \ b\sum\Sb{i,j = 1 }\\ {j > i}\endSb ^{k-1}
\pm \tilde\nabla^*\big(\check \vf(x)\vv_i,
 \check \vf(x)\vv_j\big)\wedge 
\bigwedge_{s \ne i,j}\vv_s\,.\tag 7.7
$$

Note that $\ve_{0}\wedge \check \vw(x)$ lies in the space
$\ve_{0}\wedge W^*\wedge \bigwedge^{k-2}\big((W^*)^\perp\big)$, while
$\pi(U_x\vw) - \ve_{0}\wedge \check \vw(x)$ belongs to its orthogonal
complement. Since both spaces are $D$-invariant, to prove that
$\nu\big(H(x)\vw\big) = \|\pi(DU_x\vw)\| = \|D\pi(U_x\vw)\|$ is not
less than $\rho$ for some $x\in B$ it  will suffice to show that
$\sup_{x\in B}\|D\check\vw(x)\|$ is not less than $a_0\rho$.  

Now consider the product $\ve_n\wedge\check\vw(x) $. We claim that it is enough to show that  
$$
\|\ve_n\wedge\check\vw(x) \| \ge \rho\text{  for some }x\in B\,.\tag 7.8
$$
 Indeed, since 
 $\ve_n$ is an eigenvector of $D$ with eigenvalue $a_n^{-1}$, for any
$x\in B$ the norm of $D\big(\ve_n\wedge\check\vw(x)\big)$ is not
greater than $a_n^{-1}\|D\check\vw(x)\|$. Therefore, since the
smallest eigenvalue of $D$ on $W^*\wedge \bigwedge^{k-2}\big((W^*)^\perp\big)$
is equal to $(a_*a_{n-k+1}\cdot\dots\cdot a_n)^{-1}$, the norm of
$D\check\vw(x)$ will be not less than  
$$
\split
a_n\|D\big(\ve_n\wedge\check\vw(x)\big)\| 
&\ge \frac{a_n}{a_*a_{n-k+1}\cdot\dots\cdot a_n}\|\ve_n\wedge\check\vw(x)\| \\
\un{\text{for some $x\in B$, by (7.8)}}
\ge \ \frac{\rho}{a_*a_{n-k+1}\cdot\dots\cdot a_{n-1}}\ 
\un{\text{since $a_i \ge 1$}}\ge  &\ \frac{a_0\rho}{a_0a_*a_{1}\cdot\dots\cdot a_{n-1}} \un{\text{by (5.6)}}\ge a_0\rho\,.
\endsplit
$$
Thus it remains to prove (7.8).  For this let us select the term containing $\vv_1\wedge\dots\wedge \vv_{k-2}$, and multiply (7.7) 
by $\ve_n$ as follows: 
$$
\aligned
&\ve_n\wedge\check\vw(x) =  \pm \vv^*(x)\wedge \ve_n\wedge \vv_1\wedge\dots\wedge \vv_{k-2} \\
+ \text{ other terms }&\text{where one or two of $\vv_i$, $i = 1,\dots,k-2$, are missing}\,,
 \endaligned \tag 7.9
$$
where
$$
\vv^*(x) \df \tilde\nabla^*\big(\check
\vf(x)\vv_{k-1}, a + b\check \vf(x)\vv_0 \big) = b\tilde\nabla^*\big(\check
\vf(x)\vv_{k-1}, \check \vf(x)\vv_0 \big) - \nabla^*\big(\check
\vf(x)\vv_{k-1})\,.     
$$
Because of the orthogonality of the two summands in (7.9),
and also because $\ve_n$ is orthogonal to $\vv_i$, $i = 1,\dots,k-2$, 
it follows that $\|\ve_n\wedge\check\vw(x)\|$ is not less than $\|\vv^*(x)\|$.
It follows from the first expression for $\vv^*(x)$ that $\sup_{x\in B}\|\vv^*(x)\| \ge \rho_2b$, and from the second one that  $\sup_{x\in B}\|\vv^*(x)\| \ge 
\rho_1a - 2M^2b$. An elementary computation shows that $\rho$ as defined by (7.5) is not greater than
$\min_{a^2 + b^2 \ge 1}\max(\rho_2b,\rho_1a - 2M^2b)$. This completes the 
proof of (7.8), and hence of Theorem 5.1. \qed

\heading 
{8. Completion of the proof and concluding remarks}
\endheading

 \subhead{8.1}\endsubhead First let us finish the proof of Theorem 1.1 by writing down the 

\demo{Reduction  of Theorem 1.1  to Theorems 1.3 and 1.4}
Recall that we are given  an open
subset $U$ of $\br^d$, an $n$-tuple $\vf = (f_1,\dots,f_n)$   of $C^m$ 
functions on $U$ and a function $\Psi: \bz^n\nz \mapsto(0,\infty)$ satisfying {\rm (1.1)} and(1.5).  Take  $x_0\in U$ such that $\vf$ is $l$-nondegenerate at
$x_0$ for some $l\le m$, choose $V\subset U$  as in
Theorem 1.4, and pick a 
ball $B\subset V$ containing $x_0$ such that its dilate $\tilde B$ 
(the ball with the same center as 
$B$ and twice bigger radius) is contained in $U$.
 We are going to prove that for
a.e.~$x\in B$ one has $\vf(x)\in\Cal W(\Psi)$. In other words, define $A(\vq)$ to be the set of $x\in B$ satisfying 
$
\lfloor\vf(x)\vq\rfloor < \Psi(\vq)\,;
$
we need to show that points $x$ which belong to infinitely many 
sets $A(\vq)$ form a set of measure zero.

We  proceed by induction on $n$.  If $n \ge 2$,  let us assume
that the claim is proven for any nondegenerate $(n-1)$-tuple of
functions. 
Because of the induction assumption and the fact that projections of a
nondegenerate manifold are nondegenerate, we know that  almost every $x\in B$ 
belongs to at most finitely many sets $A(\vq)$ such that $q_i = 0$ for
some $i = 1,\dots,n$. It remains to show that the same is true if one includes integer vectors
$\vq$ with 
all coordinates different from zero (if $n=1$ there is no difference, so the argument below
provides both the base and the induction step).  

Take $L$ as in (1.6a), denote by $A_\ge(\vq)$ the set of 
$x\in A(\vq)$ satisfying  (1.6d), and set $A_<(\vq)\df A(\vq)\ssm A_\ge(\vq)$. Theorem 1.3 guarantees that the measure of  $A_\ge(\vq)$ is not greater than $C_d\Psi(\vq)|B|$ whenever $\vq$ is far enough from the origin. Because of (1.1), 
the sum of measures of the sets $A_\ge(\vq)$ is finite, hence, by the   Borel-Cantelli Lemma, almost every $x\in B$ is contained in  at most finitely many sets $A_\ge(\vq)$.

Our next task is to use Theorem 1.4 to estimate the measure of the union
$$
\tsize{\bigcup_{\vq\in Q,\,2^{t_i}
\le |q_i| < 2^{t_i+1}
} } A_<(\vq)\tag 8.1
$$
for any $n$-tuple $\vt = (t_1,\dots,t_n) \in \bzp^n$ with large enough $\|\vt\| = \max_i t_i$. Observe 
that conditions (1.1) and (1.5) imply that  $\Psi(\vq) \le \left(\prod_i|q_i|\right)^{-1}$ whenever $\vq$ is far enough from the origin. It follows that $\Psi(\vq) \le 2^{-\sum_i t_i}$ whenever $\vq$ satisfies the restrictions of (8.1) with $\vt$ far enough from the origin. Therefore for such $\vt$  the set (8.1) will be contained in the set (1.7b) where one puts $\delta = 2^{-\sum_i t_i}$, $K = \sqrt{ndL}2^{\|\vt\|/2}$ and ${T}_i = 2^{{t_i+1}}$. It is straightforward to verify that inequalities (1.7a) are satisfied whenever $\|\vt\|$ is large enough; in fact one has
$$
\frac{\delta K {T}_1\cdot\dots\cdot
{T}_{n}}{\max_i{T}_i} = \frac{2^{-\sum_i t_i}\sqrt{ndL}2^{\|\vt\|/2}2^{n+\sum_i t_i}}{2^{\|\vt\|+1}} = \sqrt{ndL}2^{n-1-\|\vt\|/2}\,,
$$
which for large $\|\vt\|$ is less than $1$ but bigger than $\delta^{n+1}$. Therefore $\vre$ as in (1.7c) is equal to $\sqrt{ndL}2^{n-1}2^{-\frac1{2(n+1)}\|\vt\|}$, so, by Theorem 1.4, the measure of the set (8.1) is at most 
$$
E\big(\sqrt{ndL}2^{n-1}\big)^{-\frac1{d(2l-1)}}2^{-\frac1{2d(2l-1)(n+1)}\|\vt\|}\,.
$$
 Hence the sum of the measures of the sets (8.1) over all ${\vt\in \bzp^n}$ is finite, which implies that almost all
$x\in B$ is contained in at most finitely many such sets. To finish the proof, it remains to observe that
parallelepipeds $\{2^{t_i}
\le |q_i| < 2^{t_i+1}\}$ cover all the integer vectors $\vq$ with 
each of coordinates different from zero. 
\qed \enddemo

\subhead{8.2}\endsubhead Here is one more example of functions $\Psi$ one can consider.  For an $n$-tuple \linebreak ${\vs} = (s_1,\dots,s_n)$ with positive components, define the {\it ${\vs}$-quasinorm\/} $\|\cdot\|_{\vs}$ on $\br^n$ by\linebreak  $\|\vx\|_{\vs} \df  \max_{1\le i \le n}|x_i|^{1/s_i}$. Then, following \cite{Kl}, say that $\vy\in\br^n$ is {\sl $\vs$-$\psi$-approximable\/} if  
it belongs to $\Cal W(\Psi)$ where  
$$
\Psi(\vq) =\psi(\|\vq\|_\vs)\,.
$$ 
We will normalize $\vs$ so that $\sum_i s_i = 1$ (this way, for example, one can see that, for a non-increasing $\psi$,  any
$\vs$-$\psi$-approximable $\vy$ is 
$\psi$-MA). The choice $\vs = (1/n,\dots,1/n)$ gives the standard definition of $\psi$-approximability.  
One can also show that (1.1) holds if and only if $
\sum_{k = 1}^{\infty}{\psi(k)} < \infty$. Thus one has the following generalization of part (S) of Corollary 1.2:

\proclaim{Corollary} Let $\vf:U\to \br^n$  be as in Theorem 1.1, 
$\psi:\bn\mapsto(0,\infty)$ a
non-increasing function, and take any ${\vs} = (s_1,\dots,s_n)$ with $s_i > 0$ and $\sum_i s_i = 1$.  Then, assuming {\rm (1.2s)},   
  for        
almost all $x\in U$ the points $\vf(x)$ 
are not  $\vs$-$\psi$-approximable.
\endproclaim

\subhead{8.3}\endsubhead The idea to study the
set of points $x$ such that  $F(x) \df
\vf(x)\vq$ is close to an integer by looking at  the values of
the gradient $\nabla F(x) = \vf(x)\vq$ of $F$  has a long history. 
It was extensively used by \sp\ in his proof of Mahler's Conjecture
\cite{Sp2, Sp4}, that is, when $d = 1$ and $f_i (x) = x^i$. Also from a paper of
A.~Baker and  
W.~Schmidt \cite{BS} it follows that for some $\gamma,\vre > 0$, on a set of
positive measure the system 
$$
\cases 
|P(x)|  < \|\vq\|^{-n+\gamma} \\
 |P'(x)| < \|\vq\|^{1-\gamma- \vre} \endcases\tag 8.2
$$
(here  $P(x)$ is the polynomial  $p + q_1x + \dots + q_n
x^n$) has at most finitely many  solutions $p\in \bz$,
$\vq\in\bz^n$. This was used 
 to construct a certain regular system of real numbers and
obtain the sharp lower estimate for the \hd\ of the set 
$$
\big\{x\in\br\bigm| \lfloor q_1x + \dots + q_n
x^n\rfloor < \|\vq\|^{-\lambda}\text{ for infinitely many
}\vq\in\bz^n\big\}
$$
for $\lambda > n$. Note also that the system (8.2) is related to the
distribution of values of discriminants of integer polynomials, 
see \cite{D, Sp2, Bern}. 

In 1995 V.~Borbat \cite{Bo} proved that given any $\vre > 0$ and $0 <
\gamma < 1$, for almost all $x$ there are  at most  finitely many
solutions of (8.2). 
Now we can use Theorem 1.4 to relax the restriction $\gamma < 1$. More precisely, 
we derive the following generalization and
strengthening of the  aforementioned  result of Borbat:

\proclaim{Theorem} Let $U\subset\br^d$ be an open subset,  and
let $\vf = (f_1,\dots,f_n)$
be a nondegenerate   $n$-tuple of $C^m$ functions on $U$. Take $\vre
> 0$ and $0 < \gamma < n$. Then for
almost all $x\in V$ there exist at most finitely many solutions 
$\vq\in\bz^n$ of the system
$$
\cases 
\lfloor \vf(x)\vq \rfloor  < \Pi_{\sssize +}(\vq)^{-1+\gamma/n} \\
 \|\nabla\vf(x)\vq\| < \|\vq\|^{1-\gamma- \vre} \endcases\tag 8.3
$$
\endproclaim

\demo{Proof} As in the proof of Theorem 1.1, one can use induction to be left with integer vectors
$\vq$ with 
all coordinates different from zero. Then one estimates the measure of the union of the sets of 
solutions of (8.3) over all $\vq\in \bz^n$ with $
2^{t_i}
\le |q_i| < 2^{t_i+1}$ 
for any $n$-tuple $\vt = (t_1,\dots,t_n) \in \bzp^n$ by using Theorem 1.4 with  $\delta = 2^{(-1+\gamma/n)\sum_i t_i}$, $K = \sqrt{ndL}2^{-(1-\gamma- \vre)\|\vt\|}$ and ${T}_i = 2^{{t_i+1}}$. Inequalities (1.7a) are clearly satisfied; in particular one has
$$
\frac{\delta K {T}_1\cdot\dots\cdot
{T}_{n}}{\max_i{T}_i} = 2^{\frac\gamma n(\sum_i t_i - n\|\vt\|)}2^{-\vre\|\vt\|} \le 2^{-\vre\|\vt\|}\le 1\,.
$$
Therefore, by Theorem 1.4, the measure of the above union is at most 
$$
E \max\big(2^{\frac{-1+\gamma/n}{d(2l-1)}\sum_i t_i}, 2^{-\frac{\vre}{d(2l-1)(n+1)}\|\vt\|}\big)\,.
$$
Obviously the sums of both functions in the right hand side over all $\vt\in\bzp^n$ are finite, which completes the proof.
\qed \enddemo

The above theorem
naturally invites one to think about a possibility of  Khintchine-type 
results involving derivative estimates; that is, replacing (8.3) by,
say,  
$$
\cases 
\lfloor \vf(x)\vq \rfloor  < \Psi_1(\vq)\\
 \|\nabla\vf(x)\vq\| < \Psi_2(\vq) \endcases\tag 8.4
$$
and finding optimal conditions on $\Psi_i$ implying at most finitely many  solutions of (8.4) for almost all
$x$.

 \subhead{8.4}\endsubhead The main result of the paper (Theorem 1.1) was proved
already in the summer of 1998, but only in the case when the functions
$f_1,\dots,f_n$ are analytic. More precisely, the analytic set-up was
reduced to the case $d = 1$ (see \cite{Sp4, \S3} or 
\cite{P} for a related ``slicing'' technique). In the latter case, 
in addition to the
nondegeneracy of $\vf$, we had to assume that there exist positive
constants $C$ and $\alpha$ such that for almost all $x\in
U\subset \br$ one can find a subinterval $B$ of $U$ containing $x$ such that 
$$
\text{\rm
Span}\left(f'_1,\dots,f'_n; 
\left|\matrix f_i & f_j\\ f'_i & f_j'\endmatrix \right|,
1\le i < j \le n\right)
\text{consists of functions \cag\ on } B\,.
$$
For analytic functions $f_1,\dots,f_n$ this condition can be easily
verified by applying Corollary 3.5(a) to the basis of the above 
function space. 

In our original approach for $d = 1$ we considered sets more general
than (1.7b), 
namely  the sets
$$
\Big\{x\in B\bigm| \exists\,\vq\in\bz^n\nz\text{ such that }\cases 
\lfloor \vf(x)\vq \rfloor  < \delta \\
 |\vg(x)\vq| < K\\
|q_i| < {T}_i,\quad i = 1,\dots,n
 \endcases\ \Big\}\,,\tag 8.5
$$
with $\vf$ as in Theorem 1.4 and $\vg$ another nondegenerate
$n$-tuple of functions on $U$. We were able to prove an analogue of
Theorem 1.4 for sets (8.5) but only for $K \ge 1$, and for $n$-tuples
$\vf$ and $\vg$ 
with an additional assumption that 
$$
\text{\rm
Span}\left(g_1,\dots,g_n; 
\left|\matrix f_i & f_j\\ g_i & g_j\endmatrix \right|,
1\le i < j \le n\right)
\text{consists of  \cag\ functions}\,.\tag 8.6
$$
Instead of $U_x$ as in (5.2) we considered more general matrices 
$$U^{\vf,\vg}_x \df \left(\matrix
1 & 0 & \vf(x)  \\
0 & 1 & \vg(x) \\0 & 0 & I_n
\endmatrix \right)\,.$$. 
To prove an analogue of condition (i) of Theorem
6.2,  or, more
precisely, the statement that for some positive $C,\alpha$ and   
any subgroup $\Gamma$ of $\Lambda$ the function
$x\mapsto\|DU^{\vf,\vg}_x\Gamma\|$ is 
$(C,\alpha)$-good on some neighborhood $B$ of $x_0$, it was enough to
consider the standard basis $\{\ve_0,\ve_*\ve_1,\dots,\ve_n\}$ of
$\br^{n+1}$ and the corresponding basis 
$$
\big\{\ve_{\sssize I} \df \ve_{i_1}\wedge\dots\wedge \ve_{i_k}\mid I =
\{i_1,\dots,i_k\}\subset \{0,*,1,\dots,n\}\big\}
$$
of $\bigwedge^k(\br^{n+1})$, decompose an element $\vw$ representing
$\Gamma$ as $\vw = \sum _{\sssize I} w_{\sssize I}\ve_{\sssize I}$, write an
expansion similar to (7.3) and use (8.6). 

To prove an analogue of condition (ii), that is the statement that for any
neighborhood $B$ of  $x_0$ there exists  $\rho > 0$
 such that 
$
\sup_{x\in B}\|DU^{\vf,\vg}_x\Gamma\| \ge \rho$ for every $\Gamma \subset
\Lambda$, we used the fact that
the coefficients $w_{\sssize I}$ are integers and considered the
following two cases: 1) $w_{\sssize I} \ne 0$ for some
$I\subset \{1,\dots,n\}$, and 2) $w_{\sssize I} = 0$ for all
$I\subset \{1,\dots,n\}$. In the just described approach it was
important that $K \ge 1$. This was enough for the proof of the ``d =
1''-case of Theorem 1.1;
however, as we saw  in \S 8.3, the stronger version,
allowing arbitrarily small positive values of $K$, is important for
other applications.

\subhead{8.5}\endsubhead For completeness let us discuss the complementary divergence case of Khintchine-type theorems mentioned in the paper. It was proved by A.~Khintchine in 1924 \cite{Kh} (resp.~by A.~Groshev in 1938 \cite{G}) that a.e.~$\vy\in\br$ (resp.~$\br^n$) is $\psi$-approximable whenever $\psi$ is a non-increasing function which does not satisfy (1.1s). In 1960 W.~Schmidt \cite{S1} showed that a.e.~$\vy\in\br^n$ belongs to $\Cal W(\Psi)$ whenever the series in (1.1) diverges (note that there are no monotonicity restrictions on $\Psi$ unless  $n = 1$). It seems plausible to conjecture the divergence counterpart of Theorem 1.1, namely that for $\vf:U\mapsto \br^n$ as in Theorem 1.1 and $\Psi$ satisfying (1.5) but not (1.1), the set $\{x\in U\mid \vf(x)\in\Cal W(\Psi)\}$ has full measure. For functions $\Psi$ of the form (1.2s) this can be done using Theorem 1.4 and the method of {\sl regular
systems}, which dates back to \cite{BS} and has been extensively used 
in the existing proofs of divergence Khintchine-type results for 
special classes of manifolds \cite{DRV2, DRV3, BBDD, Be1, Be3, Be4}. 


\heading{Acknowledgements}
\endheading 

A substantial part of this work was done during the authors' stays at
the University of Bielefeld in 1998 and 1999. These stays were
supported by SFB-343 and Humboldt Foundation. The paper was completed during the Spring 2000 programme on Ergodic Theory, Geometric Rigidity and Number Theory at the Isaac Newton Institute of Mathematical Sciences. Thanks are also due to
V.~Beresnevich and M.~Dodson for useful remarks.

\enddocument

\end

\enddocument